\documentclass[10pt,onecolumn,journal,draftclsnofoot]{IEEEtran}

\DeclareSymbolFont{bbold}{U}{bbold}{m}{n}
\DeclareSymbolFontAlphabet{\mathbbold}{bbold}

\usepackage{cleveref}
\usepackage{mathtools}

\usepackage{graphicx,color,url,sidecap}
\usepackage{amsmath,amssymb,mathrsfs,amsfonts}
\usepackage{enumitem}

\usepackage{hyperref}

\usepackage{doi}

\newtheorem{theorem}{Theorem}[section]
\newtheorem{lemma}[theorem]{Lemma}

\newtheorem{definition}{Definition}[section]
\newtheorem{remark}{Remark}

\newcommand{\until}[1]{\{1,\dots, #1\}}
\newcommand{\fromto}[2]{\{#1,\dots, #2\}}

\newcommand{\setdef}[2]{\{#1 \; | \; #2\}}

\newcommand\oprocendsymbol{\hbox{$\square$}}
\newcommand\oprocend{\relax\ifmmode\else\unskip\hfill\fi\oprocendsymbol}

\newcommand{\vect}[1]{\mathbbold{#1}}
\newcommand{\vectorones}[1][]{\vect{1}_{#1}}
\newcommand{\vectorzeros}[1][]{\vect{0}_{#1}}

\newcommand{\ml}{\left[ \begin{array}{cccccccccccccccccccccccccc}}
\newcommand{\mr}{\end{array} \right]}

\graphicspath{{./fig/}}

\newcommand{\real}{\mathbb{R}}
\renewcommand{\natural}{\mathbb{N}}

\newcommand{\diag}{\operatorname{diag}}
\newcommand{\interior}{\operatorname{interior}}
\newcommand{\simplex}[1]{\Delta_{#1}}

\newcommand{\argmin}{\ensuremath{\operatorname{argmin}}}

\title{Opinion Dynamics and Social Power Evolution:\\ A
  Single-Timescale Model}
\author{Peng Jia, Noah E.\ Friedkin, Francesco
  Bullo, \IEEEmembership{Fellow, IEEE}
  \thanks{This material is based upon work
    supported by the U. S.\ Army Research Laboratory and the
    U. S.\ Army Research Office under grants number W911NF-15-1-0274
    and W911NF-15-1-0577.  The content of the information does not
    necessarily reflect the position or the policy of the Government,
    and no official endorsement should be inferred.}
  \thanks{Peng Jia, Noah Friedkin, and Francesco Bullo are
    affiliated with the Center for Control, Dynamical Systems and
    Computation, University of California at Santa Barbara, 
    {\tt pengjwt@gmail.com}, 
    {\tt bullo@engineering.ucsb.edu}, 
    {\tt friedkin@soc.ucsb.edu}.}}

\begin{document}
\maketitle

%\IEEEcompsoctitleabstractindextext{%

\begin{abstract}
  This paper studies the evolution of self-appraisal and social power,
  for a group of individuals who discuss and form opinions. We
  consider a modification of the recently proposed DeGroot-Friedkin
  (DF) model, in which the opinion formation process takes place on
  the same timescale as the reflected appraisal process; we call this
  new model the single-timescale DF model. We provide a comprehensive
  analysis of the equilibria and convergence properties of the model
  for the settings of irreducible and reducible influence networks.
  For the setting of irreducible influence networks, the
  single-timescale DF model has the same behavior as the original
  DF model, that is, it predicts among other things that the social
  power ranking among individuals is asymptotically equal to their
  centrality ranking, that social power tends to accumulate at the top
  of the centrality ranking hierarchy, and that an autocratic (resp.,
  democratic) power structure arises when the centrality scores are
  maximally nonuniform (resp., uniform).
  For the setting of reducible influence networks, the
  single-timescale DF model behaves differently from the original DF
  model in two ways.
  First, an individual, who corresponds to a reducible node in a reducible
  influence network, can keep all social power in the single-timescale DF
  model if the initial condition does so, whereas its social power 
  asymptotically vanishes in the original DF model.
  Second, when the associated network has multiple sinks, the two 
  models behave very differently: the original DF model has a single
  globally-attractive equilibrium, whereas any partition of social
  power among the sinks is allowable at equilibrium in the
  single-timescale DF model.
\end{abstract}


%}

%\IEEEdisplaynotcompsoctitleabstractindextext

\IEEEpeerreviewmaketitle

%------------------------------------------
\section{Introduction}\label{sec:introduction}
%------------------------------------------

\paragraph*{Problem description and motivation}
This article focuses on a model for the evolution of social power and
self-appraisal in an influence network.  The model combines an opinion
dynamics process from network systems and a reflected appraisal
process from applied psychology. The model is a variation of a
recently-proposed dynamical system, called the DeGroot-Friedkin (DF)
model, proposed and characterized in~\cite{PJ-AM-NEF-FB:13d}; in this
proposed variation, the opinion dynamics process takes place on the
same timescale as the reflected appraisal process.  In other words, in
the (original) DF model reflected appraisal played out over an issue discussion 
sequence where opinion consensus was reached on each issue; here it 
plays out during the opinion influence process on a single issue.
The model we study in this paper was also independently proposed and
studied by Xu \emph{et al.}~\cite{ZX-JL-TB:15}.  The purpose of this
article is to provide a rigorous and comprehensive analysis of the
asymptotic behavior of the proposed model and to compare it with the
DF model.

\paragraph*{Literature review}
Influence networks and opinion formation processes have been the
subject of a rich literature, starting with the averaging model
proposed by French in \cite{JRPF:56}, studied also by Harary in
\cite{FH:59} and DeGroot in \cite{MHDG:74}, the Abelson
model~\cite{RPA:64}, the Friedkin-Johnsen model
\cite{NEF-ECJ:99,NEF-AVP-RT-SEP:16}, and the Hegselmann-Krause model
\cite{RH-UK:02} among others. Empirical evidence in support of the
averaging model (including its variations) is described
in~\cite{NEF-ECJ:11,AC-HL-JX:16}.  These models are now standard in
surveys and textbooks such
as~\cite{CC-SF-VL:09,MOJ:10,AVP-RT:17,FB:17}.

Recently, by combining the DeGroot model of opinion dynamics and a
reflected appraisal mechanism, Jia \emph{et al.}
\cite{PJ-AM-NEF-FB:13d} proposed a DF model (the DF
model) to describe the evolution of individuals' self-appraisal and
social power in a network along an issue sequence. Empirical evidence
in support of the reflected appraisal mechanism and other aspects of
the DF model is provided in~\cite{NEF-PJ-FB:14n}, which present a
remarkable suite of issue-sequence effects on influence network
structure consistent with theoretical predictions.

Building on the modeling ideas in~\cite{PJ-AM-NEF-FB:13d}, several
extensions and variations have been proposed recently.  For example,
Mirtabatabaei \emph{et al.} extended the DF model to include stubborn
agents who have attachment to their initial opinions in
\cite{AM-PJ-NEF-FB:13x}.  A continuous-time self-appraisal model was
introduced by Chen \emph{et al.} in \cite{XC-JL-MAB-ZX-TB:17}.
Considering time-varying doubly stochastic influence matrices, Xia
\emph{et al.}  \cite{WX-JL-KHJ-TB:16} investigated the convergence
rate of the modified DF model, which was proven to converge
exponentially fast.
Very recent submissions (essentially simultaneous with and independent
of this article)
include~\cite{MY-JL-BDOA-CY-TB:17,GC-XD-NEF-FB:16j,ZA-RF-AH-YC-TTG:17};
specifically, the works~\cite{MY-JL-BDOA-CY-TB:17,GC-XD-NEF-FB:16j}
deal with time-varying (deterministic or stochastic) influence
networks and the article~\cite{ZA-RF-AH-YC-TTG:17} provides novel
stability analysis methods for nonlinear Markov chains (motivated by
the DF model).

Finally and notably, motivated by~\cite{PJ-AM-NEF-FB:13d},
Xu \emph{et al.} \cite{ZX-JL-TB:15} proposed a
\emph{modified DF model} where the social power is updated without
waiting for opinion consensus on each issue to take place, i.e., the
local estimation of social power is truncated. In this sense, the
time-constant of the opinion dynamic process is now the same as that
for the reflected appraisal process.
The analysis of the equilibrium points and their attractivity
properties was given in \cite{ZX-JL-TB:15} only for the setting where
the interaction matrix is doubly stochastic.  This is the model
studied in this paper under the name ``single-timescale
DeGroot-Friedkin (DF) model.''

\paragraph*{Statement of contributions}
Section~\ref{sec:iF-model} introduces the main modeling assumptions
and the definition of the single-timescale DF model.
Section~\ref{sec:dynamics-irreducible} provides a comprehensive
analysis of the proposed model for irreducible influence
networks. Specifically, Theorem~\ref{thm:row-stochastic-star}
characterizes the system behavior over influence networks with star
topology and Theorem~\ref{thm:row-stochastic-general} treats the
general case.  The latter theorem subsumes the specific setting of
doubly-stochastic influence networks. Lemma~\ref{lemma:more-features}
characterizes the relationship with the DF model: the two models (over
irreducible influence networks) converge to the same equilibria and
therefore predict the same phenomena, e.g., social power ranking equal
to an appropriate centrality ranking and social power accumulation at
the top.
Next, Section~\ref{sec:dynamics-reducible} treats the setting of
reducible influence networks. Theorem~\ref{thm:reducible-reachable}
shows that the single-timescale DF model behaviors similarly to the DF
model over reducible influence networks with globally reachable nodes,
but its set of equilibrium points contain all vertices of a simplex,
including the cases that reducible nodes have all social power. In
contrast, the reducible nodes loss their social power asymptotically
in the DF model. Theorem~\ref{thm:reducible-no-reachable}
considers the most general case where the associated network has
multiple sinks and the two models behave very differently: the DF 
model has a single globally-attractive equilibrium,
whereas any partition of social power among the sinks is allowable at
equilibrium in the single-timescale DF model.
Finally, Section~\ref{sec:conclusion} contains some final remarks and 
all proofs are in the Appendices in the supplementary file.

In summary, we believe that these results are meaningful as they
extend the validity and scope of the original analysis. It is
important to establish the weakest possible conditions under which
social power and self-appraisal evolve in a way comparable (or
identical) to that predicted by the DF model.  This paper, together
with other efforts on time-varying influence networks, establishes
some robustness in the dynamic behavior with respect to modeling
uncertainties.

%------------------------------------------
\section{The single-timescale DF model}
\label{sec:iF-model}
%------------------------------------------

In this section we introduce and motivate the dynamical model for the
evolution of the social influence network where social opinions and
social power evolve simultaneously. This model combines the concepts
of the DeGroot model for the dynamics of opinions over a single issue
and of the Friedkin model for the dynamics of self-weight and social
power over a sequence of issues.

We consider a group of $n\geq2$ individuals who discuss an
issue according to a DeGroot opinion formation model with an
influence matrix $W$.  Assume that individual opinions about the issue are 
described by a trajectory $t\mapsto y(t)\in\real^n$ that is determined by 
the DeGroot averaging model
\begin{equation}
  \label{eq:DeGroot-simple}
  y(t+1) = W y(t),  \quad t=0, 1, 2,\dots,
\end{equation}
with given initial conditions $y_i(0)$ for each individual $i$. Here,
the influence matrix $W$ is row-stochastic, i.e., each entry of $W$ is
non-negative and each row sum of $W$ equals $1$. By~\eqref{eq:DeGroot-simple}, 
each individual $i$ updates its opinion according to the convex combination:
\begin{equation*}
  \label{eq:DeGroot-simple-components}
  y_i(t+1) = w_{ii} y_{i}(t) +
  \sum\nolimits_{j=1,j\neq{i}}^n w_{ij}y_{j}(t).
\end{equation*}
From a psychological viewpoint, the diagonal and the off-diagonal
entries of an influence matrix $W$ play conceptually distinct roles.
Specifically, the diagonal \emph{self-weight} $w_{ii}$ is the
individual's self-appraisal (e.g., self-confidence, self-esteem,
self-worth) and corresponds to the extent of closure to interpersonal
influence of the $i$th individual. Instead, the off-diagonal entries
$w_{ij}$, $j\neq{i}$, are \emph{interpersonal weights} that the $i$th
individual \emph{accords} to other individuals.

For simplicity of notation, we adopt the shorthand $x_i\in[0,1]$ to denote 
the self-weight $w_{ii}$ of the $i$th individual. Because $1-x_i$ is the 
aggregated influence on the $i$th individual of all other individuals, we 
may decompose the off-diagonal entries as $w_{ij} = (1-x_i) c_{ij}$, where 
the coefficients $c_{ij}$ are the \emph{relative interpersonal weights} that 
the $i$th individual accords to other individuals. Given $c_{ii}=0$, the 
matrix $C$, called the \emph{relative interaction matrix} is row-stochastic 
with zero diagonal. Our construction assumes that the matrix $C$ is constant.
With these notations and assumptions, a time-dependent influence matrix
is written as
\begin{equation}
  W(x(t))  = \diag(x(t))+(I_n-\diag(x(t)))C,
  \label{def:decomposition}
\end{equation}
and the opinion dynamic process~\eqref{eq:DeGroot-simple} is
rewritten as
\begin{equation*}
  \label{eq:DeGroot}
  y(t+1) = W(x(t)) y(t), \quad t=0, 1, 2,\dots.
\end{equation*}
If $C$ is further assumed to be irreducible, the Perron-Frobenius 
Theorem for non-negative matrices implies that the influence matrix $W(x)$ 
with $x\geq 0$ admits a unique left eigenvector $w(x)^\top\geq0$ associated with 
the eigenvalue $1$, with non-negative entries. We may normalize $w(x)$ so that
$w(x)\in\simplex{n}$. We refer to this row vector $w(x)^\top$ the 
\emph{dominant left eigenvector} of $W(x)$. If $W(x)$ is aperiodic additionally, 
then
\begin{equation*}
  \lim_{t\to\infty}W(x)^t = \vectorones[n]w(x)^\top.
\end{equation*}

Our model is completed by formulating how the self-weights
$t\mapsto x(t)$ evolve during the opinion formation. By adopting to the 
psychological concept of \emph{reflected appraisal}, we assume that 
individual social powers are adjusted along group discussions and the 
self-weight of an individual is set equal to the social power that the 
individual exercised over the influence network. We proposed a natural 
dynamical process~\cite{PJ-AM-NEF-FB:13d} that allows each individual to 
accurately estimate her perceived power. The dynamical process is distributed 
in the sense that each individual only needs to interact with her influenced 
neighbors (i.e., those who accord positive interpersonal weights to the 
individual). By assuming that she is aware of the direct interpersonal 
weights accorded to her and the perceived powers of her influenced 
neighbors, each individual updates her perceived power as a convex 
combination of her own and her influenced neighbors' perceived powers. 
That is, in each discussion iteration, each individual $i$ estimates 
her perceived power $p_i(t)$ according to
  \begin{equation}
  \label{eq:each-individual-power-evolution}
  \begin{split}
    p_i(t+1) &= w_{ii}(t) p_{i}(t) +\sum\nolimits_{j=1,j\neq{i}}^n
    w_{ji}(t)p_{j}(t),  \\
    &\hspace{3.5cm} t=0, 1, 2,\dots,
   \end{split}
 \end{equation}
or, equivalently, $p(t+1)=W(t)^\top p(t)$, where $W(t)$ represents the 
influence matrix associated to the issue discussion process. 
By assuming the self-weight of an individual is set equal to the social 
power that the individual exercised over the influence network, we have
$p(t)=x(t)$ for all $t$. In short, the appraisal update mechanism 
``self-weight $:=$ relative control from the influence network'' is 
written as
\begin{equation} 
\label{eq:model-1}
  x(t+1) = W(x(t))^\top x(t), \quad t=0, 1, 2,\dots.
\end{equation}
Because of the row stochastic $W(t)$, the sum of all elements
of $x(t)$ is constant. Therefore, it is convenient to assume that the
self-weight vector $x(t)$ takes value in $\simplex{n}$ for all time $t$.

Given a vector $x=[x_1, \dots, x_n]$, we denote $x^2=[x^2_1, \dots,
  x^2_n]$ with a slight abuse of notation and then $\vect{e}^2_i=\vect{e}_i$. We conclude this modeling
discussion with a summary definition.
\begin{definition}[The single-timescale DF model for the evolution of 
social influence networks]
  \label{def:IE-model}
  Consider a group of $n\geq2$ individuals discussing an issue. 
  Let a row-stochastic zero-diagonal irreducible matrix
  $C$ be the relative interaction matrix encoding the relative
  interpersonal weights among the individuals. The single-timescale DF 
  model for the evolution of the self-weights $t\mapsto x(t)\in\simplex{n}$ 
  is defined as
  \begin{equation}
    \label{eq:sys_IF}
    \begin{split}
      x(t+1)=F(x(t)) &:= C^\top x(t) + (I-C^\top)x^2(t)\\ &\phantom{:}= C^\top \left(x(t)-x^2(t)\right)+x^2(t).
    \end{split}
  \end{equation}
\end{definition}

In this paper, we aim to (i) characterize the existence, stability, and
region of attraction of the equilibria for the single-timescale DF
model, and (ii) compare the behavior of the single-timescale DF model 
with the DF model. Based upon
Definition~\ref{def:IE-model} of the single-timescale DF model and the
definition of the DF model in~\cite{PJ-AM-NEF-FB:13d}, both models try
to describe and predict evolving social-power configures within a
social network and try to explain when and why specific configures of
self-weights (e.g., $x=\vect{e}_i$, namely autocratic configuration,
or $x=\frac{1}{n}\vectorones[n]$, namely democratic configuration) are
attractive. Nevertheless, the evolution of the single-timescale DF model
is defined on a single issue discussion, that is, the process of
opinion dynamics and the process of reflected appraisal take place
over comparable timescales (in sense that the individual self-weight
$x_i$ is set equal to the individual perceived power $p_i$
in~\eqref{eq:each-individual-power-evolution} right after each opinion
discussion iteration). Compared with that, the DF model is applied to
group discussion on a sequence of issues, that is, the timescales for
the two processes are separate: the opinion dynamics are faster than
the reflected appraisal dynamics in the influence network. In other
words, opinion consensus is achieved before individual self-weights
are updated.

%------------------------------------------
\section{The single-timescale DF model over irreducible influence networks}
\label{sec:dynamics-irreducible}
%------------------------------------------

In this section we begin the mathematical analysis of the single-timescale DF model. 
We consider two meaningful situations where the relative interaction matrix 
$C$ has star topology and where the digraph associated to $C$ is row-stochastic 
(including its special case where $C$ is doubly-stochastic). We will show that the 
first situation leads to the emergence of an autocratic power structure with a 
single leader from all initial conditions, and the second situation leads to the 
general convergence of self-weight configures, including the emergence of a 
democratic power structure for doubly-stochastic $C$.

\subsection{Interactions with star topology and autocratic influence networks}

Consider the first case where the digraph associated to
the relative interaction matrix has star topology. We assume $n\geq3$
because the case $n=2$ is trivial (where $C$ is necessarily symmetric and
doubly-stochastic).

\begin{theorem}[Single-timescale DF model with star topology]
  \label{thm:row-stochastic-star} 
  For $n\geq3$, consider the single-timescale DF dynamical system
  $x(t+1)=F(x(t))$ defined by a relative interaction matrix
  $C\in\real^{n\times{n}}$ that is row-stochastic, irreducible, and 
  has zero diagonal.  
  If $C$ has star topology with center node $1$, then
  \begin{enumerate}
    
  \item[\rm{(i)}]
    (\textbf{Equilibria:}) \label{fact:equilibria-rs-star} the fixed
    points of $F$ are the autocratic vertices
    $\{\vect{e}_1,\dots,\vect{e}_n\}$, and

  \item[\rm{(ii)}] (\textbf{Convergence property:}) 
    \label{fact:stability-rs-star} 
    for all non-autocratic initial conditions 
    $x(0)\in\simplex{n}\setminus\{\vect{e}_1,\dots,\vect{e}_n\}$, the
    self-weights $x(t)$ converges asymptotically to the
    autocratic configuration $\vect{e}_1$ as $t\to\infty$.
\end{enumerate}
\end{theorem}

The result of Theorem~\ref{thm:row-stochastic-star} can be interpreted
as follows.  For the single-timescale DF model associated with star
topology, the autocrat is predicted to appear on the center node along
the opinion formation process -- independently of the initial values
in almost all scenarios (except those autocratic states corresponding
to the equilibrium points of the system~\eqref{eq:sys_IF}). This is
identical to the DF model.

\subsection{Row-stochastic interactions and democratic influence networks} 

Now we consider the second case where the relative interaction matrix $C$ 
is row-stochastic. Note that $C=\begin{bmatrix} 0 & 1\\ 1& 0\end{bmatrix}$ 
for $n=2$ is such that, for any $(x_1,x_2)\in\simplex{2}$ with strictly positive 
components, $F$ in~\eqref{eq:sys_IF} always satisfies $F(x_1,x_2)=(x_1,x_2)$. 
We therefore discard this trivial case $n=2$.

\begin{theorem}[Single-timescale DF model with row-stochastic interactions]
\label{thm:row-stochastic-general} 
  For $n\geq3$, consider the single-timescale DF dynamical system
  $x(t+1)=W(x(t))^\top x(t)$ defined by a relative interaction matrix
  $C\in\real^{n\times{n}}$ that is row-stochastic, irreducible, and has
  zero diagonal.
  Assume that the digraph $G(C)$ associated to $C$ does not have star
  topology and let $c^\top$ be the dominant left eigenvector of $C$.  Then
  \begin{enumerate}
  \item[\rm{(i)}]
    (\textbf{Equilibria:}) \label{fact:equilibria-on-simplex} the set
    of fixed points of $F$ is $\{\vect{e}_1,\dots,\vect{e}_n, x^*\}$,
    where $x^*$ lies in the interior of the simplex $\simplex{n}$ and
    the ordering of the entries of $x^*$ is equal to the ordering of
    the entries of $c$, and

    \item[\rm{(ii)}] (\textbf{Convergence
      property:}) \label{fact:general-stability} for all
      non-autocratic initial conditions
      $x(0)\in\simplex{n}\setminus\{\vect{e}_1,\dots,\vect{e}_n\}$,
      the self-weights $x(t)$ exponentially converges to the equilibrium
      configuration $x^*$ as $t\to\infty$.
\end{enumerate}
\end{theorem}

Based upon the proof of Theorem~\ref{thm:row-stochastic-star} (i) 
in Appendix~\ref{pf:thm:row-stochastic-star} and the proof of 
Theorem~\ref{thm:row-stochastic-general} (i) 
in Appendix~\ref{pf:thm:row-stochastic-general} , we immediately have the 
following extended results.
\begin{lemma}[Relationship with the DF model over irreducible networks]\label{lemma:more-features}
  Given the same $C$ and the same non-autocratic initial state $x(0)$, the
  dynamical system~\eqref{eq:sys_IF} for the single-timescale DF model
  converges to the same equilibrium as the dynamical system for the DF
  model in~\cite{PJ-AM-NEF-FB:13d}. 
%  Consequentially, the ordering of
%  the vector components of the non-autocratic equilibrium $x^*$ of the
%  system~\eqref{eq:sys_IF} with a non-star topology is consistent with
%  that of $c$. 
  Consequentially, the social power in the
  dynamical system~\eqref{eq:sys_IF} is accumulated to the
  individuals $\{i\}$ in the social network with high $\{c_i\}$
  values.
\end{lemma}

The social power accumulation statement of
Lemma~\ref{lemma:more-features} is directly from the same property of
the DF system. (See details in~\cite{PJ-AM-NEF-FB:13d}.)

Although the opinion formulation and social power evolution timescales
for the single-timescale DF model and the DF model are different, the
equilibrium results of Theorem~\ref{thm:row-stochastic-general} and
Lemma~\ref{lemma:more-features} are identical to those of the DF model
with an identical $C$: the equilibrium properties from both models are
uniquely determined by the dominant left eigenvector $c^\top$ of $C$ (where $c$
can be called \emph{eigenvector centrality scores} as
from~\cite{PJ-AM-NEF-FB:13d}). In details, given an irreducible $C$
without star topology, the vector of self-weights $x(s)$ in the
single-timescale DF model converges to a unique equilibrium value
$x^*$ for all initial conditions, except the autocratic states. This
equilibrium value $x^*$ is uniquely determined by the eigenvector
centrality score $c$. The entries of $x^*$ are strictly positive and
have the same ordering as that of $c$, that is, if the centrality
scores satisfy $c_i>c_j$, then the equilibrium social power $x^*$
satisfies $x_i^*>x_j^*$, and if $c_i=c_j$, then $x_i^*=x_j^*$. The
model exhibits an interesting phenomenon similarly as from the DF
model: an accumulation of social power in the central nodes of the
network. The accumulation phenomenon is most evident for the star
topology case: the center individual with $c_i=0.5$ has a self-weight
of $1$, and all other individuals have $0$ social powers even they may
have strictly positive centrality scores. In contrast, if $C$ is
doubly-stochastic, Theorem~\ref{thm:row-stochastic-general} and
Lemma~\ref{lemma:more-features} imply the self-weights of the
single-timescale DF system exponentially converge to a democratic
configure where the social power of each individual is uniform.

\subsection*{Numerical examples on irreducible networks}
\label{subsec:numerical}
In this section we compare the dynamical behavior of the single-timescale DF 
model~\eqref{eq:sys_IF} with that of the DF model in~\cite{PJ-AM-NEF-FB:13d} 
over an influence network with star topology and over a general irreducible 
influence network.

\paragraph{A network with star topology}
We first simulate the self-weight evolution in a network with star topology $C$. 
\begin{equation}
   C=\left[\begin{matrix}0 &\frac{1}{9} &\frac{1}{9} &\frac{1}{9} &\frac{1}{9} &\frac{1}{9}&\frac{1}{9}& \frac{1}{9} &\frac{1}{9}& \frac{1}{9}\\
   1 &0& 0 &0& 0 &0  &0 &0& 0 &0\\
   1 &0& 0 &0& 0 &0  &0 &0& 0 &0\\
   1 &0& 0 &0& 0 &0  &0 &0& 0 &0\\
   1 &0& 0 &0& 0 &0  &0 &0& 0 &0\\
   1 &0& 0 &0& 0 &0  &0 &0& 0 &0\\
   1 &0& 0 &0& 0 &0  &0 &0& 0 &0\\
   1 &0& 0 &0& 0 &0  &0 &0& 0 &0\\
   1 &0& 0 &0& 0 &0  &0 &0& 0 &0\\
   1 &0& 0 &0& 0 &0  &0 &0& 0 &0  
   \end{matrix}\right].
   \label{eq:C-star}
\end{equation}

\begin{figure}[htp]
\begin{center}
    \includegraphics[width=.8\textwidth]{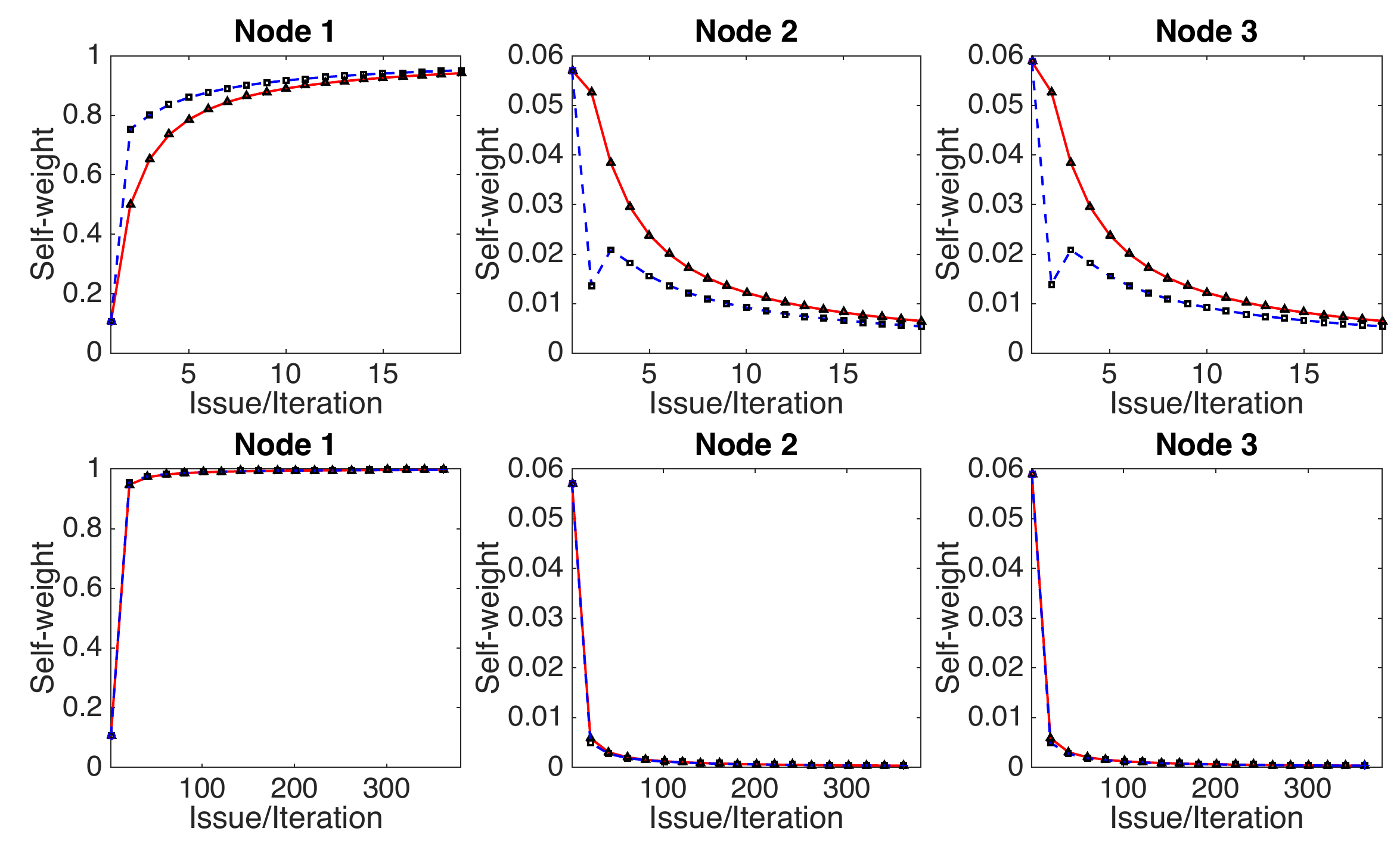} 
    \caption{Self-weight evolution for a network with star topology: we simulate 
    both dynamics of the single-timescale DF model and of the DF model 
    with the same initial conditions; we display the trajectories of 3 nodes. The 
    dot lines are related to the single-timescale DF model and the solid lines are 
    related to the DF model. The top figures show the short-term 
    behaviors and the bottom figures show the long-term dynamics.}
    \label{fig:str-t12}
    \end{center}
\end{figure}

In a network associated with $C$ given in~\eqref{eq:C-star}, the
dynamical trajectories of the self-weights generated by the
single-timescale DF model and by the DF model are illustrated in 
Figure~\ref{fig:str-t12}. Two models converge
to the same equilibrium. Specifically, individual $1$ has $1/2$ eigenvector 
centrality score and her equilibrium self-weight (social power) is $1$; the rest 
$9$ individuals have $1/18$ eigenvector centrality score for each and their 
equilibrium self-weights (social powers) are $0$. The social power accumulation 
phenomenon is most evident in such a network with star topology.
   
\begin{figure}[htp]
%  \begin{minipage}[t]{0.35\linewidth}
    \centering
     \includegraphics[width=.4\textwidth]{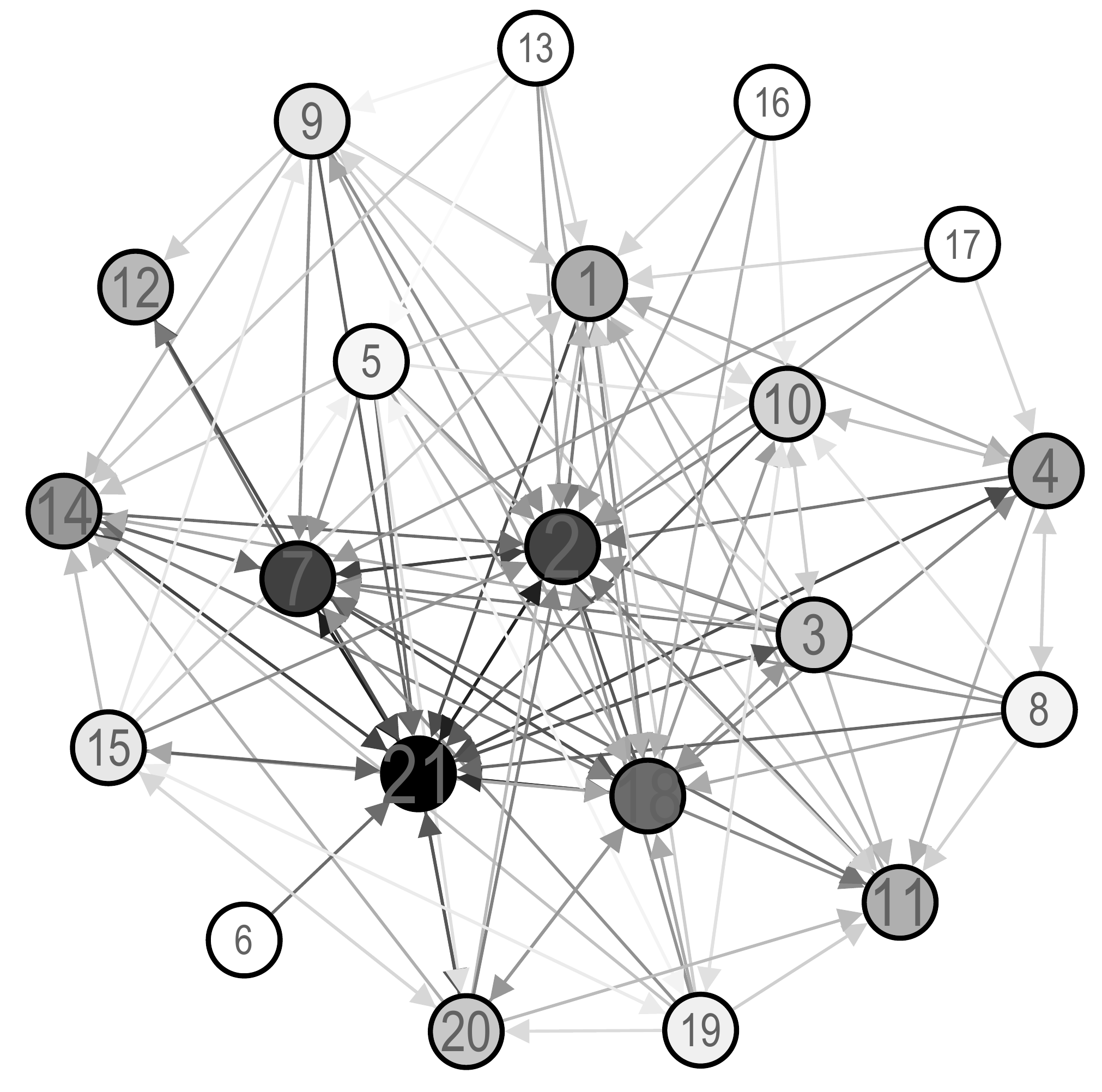}
    \caption{Krackhardt's advice network with all $21$ nodes. 
    The color gradation of the nodes and the font size of the 
    node labels represent $c_i$.}
    \label{fig:KAN-topology}
%  \end{minipage}
 % \hspace{.4cm}
%  \begin{minipage}[t]{0.47\linewidth}
\end{figure}

\paragraph{Reduced Krackhardt's advice network}
Krackhardt's advice network, as illustrated in Figure~\ref{fig:KAN-topology}, is 
based upon a US manufacturing organization, which represents $21$ managers 
and a directed advice network $C$ characterizing who sought advice from 
whom~\cite{DK:87}. If individual $i$ asks for advice from $n_i$ different 
individuals, then we assume that $c_{ij}=1/n_i$ for $j$ in these $n_i$ individuals, 
and $c_{ik}=0$ for all other individuals $k$. (See a similar example 
in~\cite{MOJ:10}.) Moreover, self-weighting is not considered in $C$, that
is, $c_{ii}=0$ for all $i\in \until{21}$.

\begin{figure}[htp]
%  \begin{minipage}[t]{0.47\linewidth}
    \centering
     \includegraphics[width=.4\textwidth]{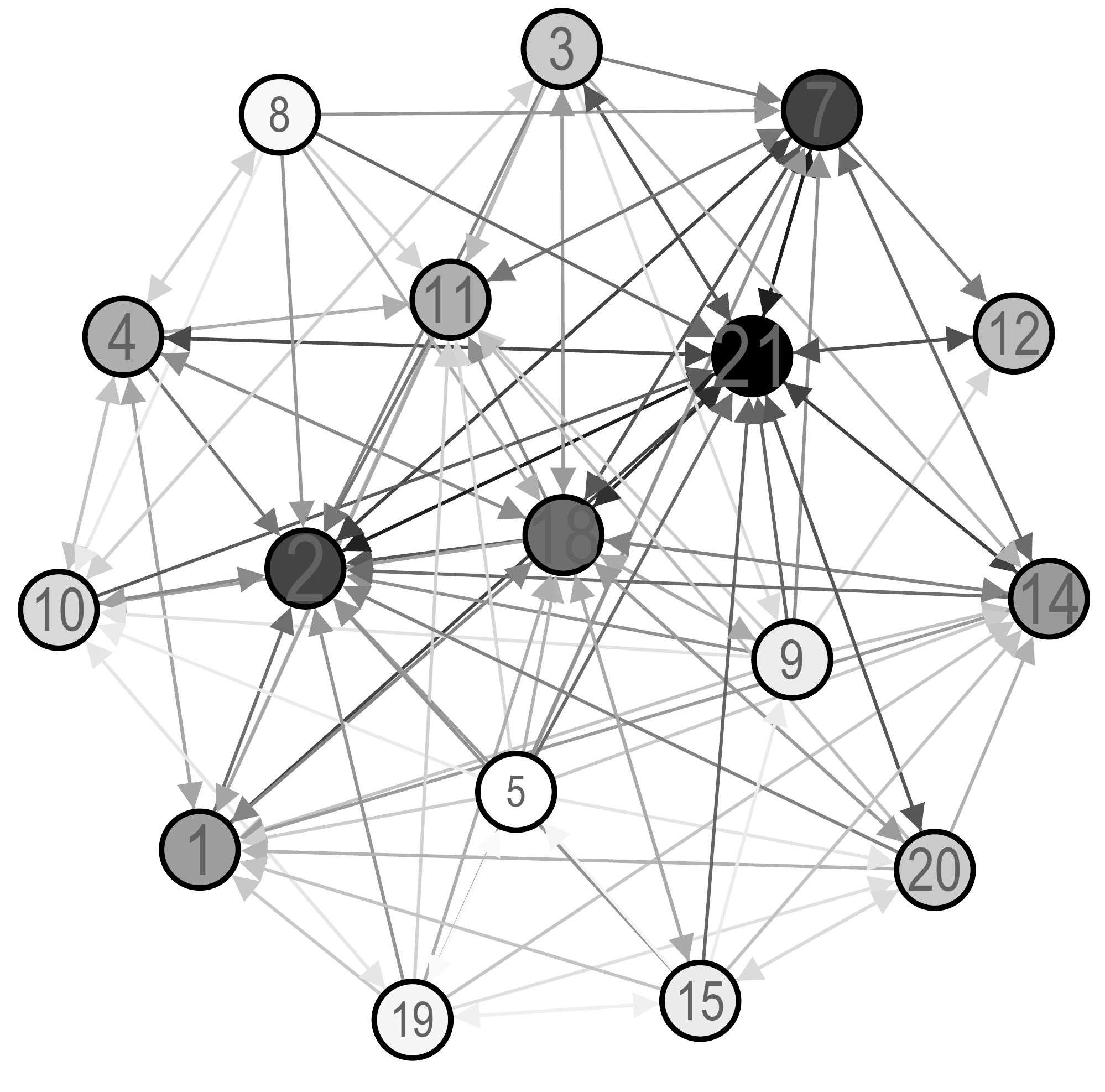}
    \caption{Reduced Krackhardt's advice network with $17$ nodes: the
      nodes $6$, $13$, $16$ and $17$ from the complete Krackhardt's advice
      network as in Figure~\ref{fig:KAN-topology} are excluded. The 
      color gradation of the nodes and the font size of the node labels 
      represent $c_i$.}
    \label{fig:red-KAN-topology}
    \end{figure}

\begin{figure}[htp]
    \centering
    \includegraphics[width=.6\textwidth]{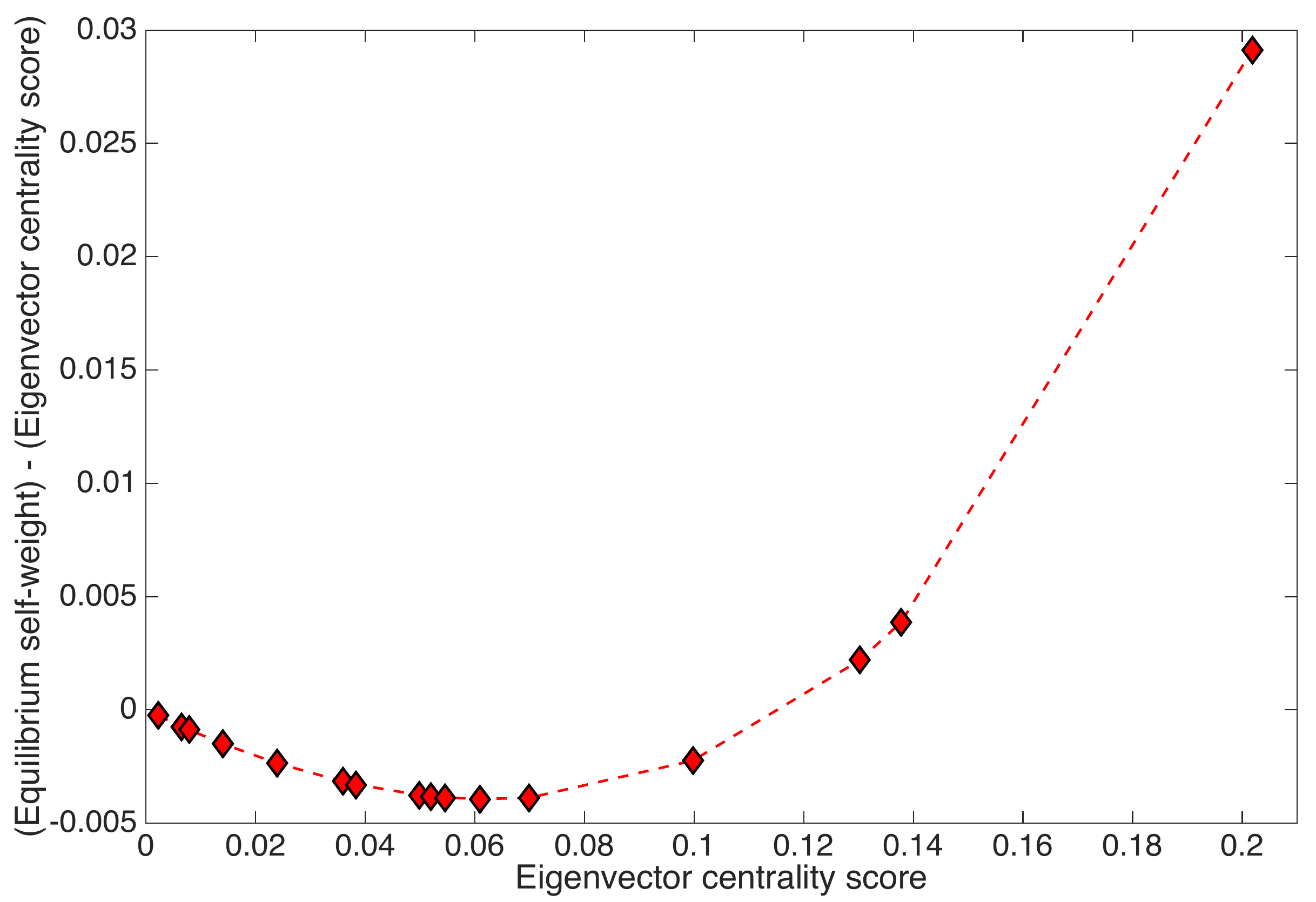}
    \caption{Comparison between the eigenvector centrality scores and
      the equilibrium self-weights for the reduced Krackhardt's advice 
      network: the social power accumulation.}
    \label{fig:red-KAN-cvsx}
%  \end{minipage}
\end{figure}

The complete Krackhardt's network includes four managers (i.e., individuals 
$6$, $13$, $16$ and $17$) from whom no other individual requests advice. 
Hence, the complete Krackhardt's network is reducible. Here, we simulate the
single-timescale DF model on a 
reduced Krackhardt's advice network (as shown in 
Figure~\ref{fig:red-KAN-topology}) without these four nodes. 
The social power accumulation phenomenon within the reduced Krackhardt's 
advice network is demonstrated in Figure~\ref{fig:red-KAN-cvsx}.
We may also check from the simulation that the 
ordering of the vector components of $x^*$ is consistent with that of $c$, that is, 
$x_i^*> x_j^*$ if and only if $c_i> c_j$ for $i, j \in \until{17}$.

\begin{figure}[htp]
\begin{center}
    \includegraphics[width=.8\textwidth]{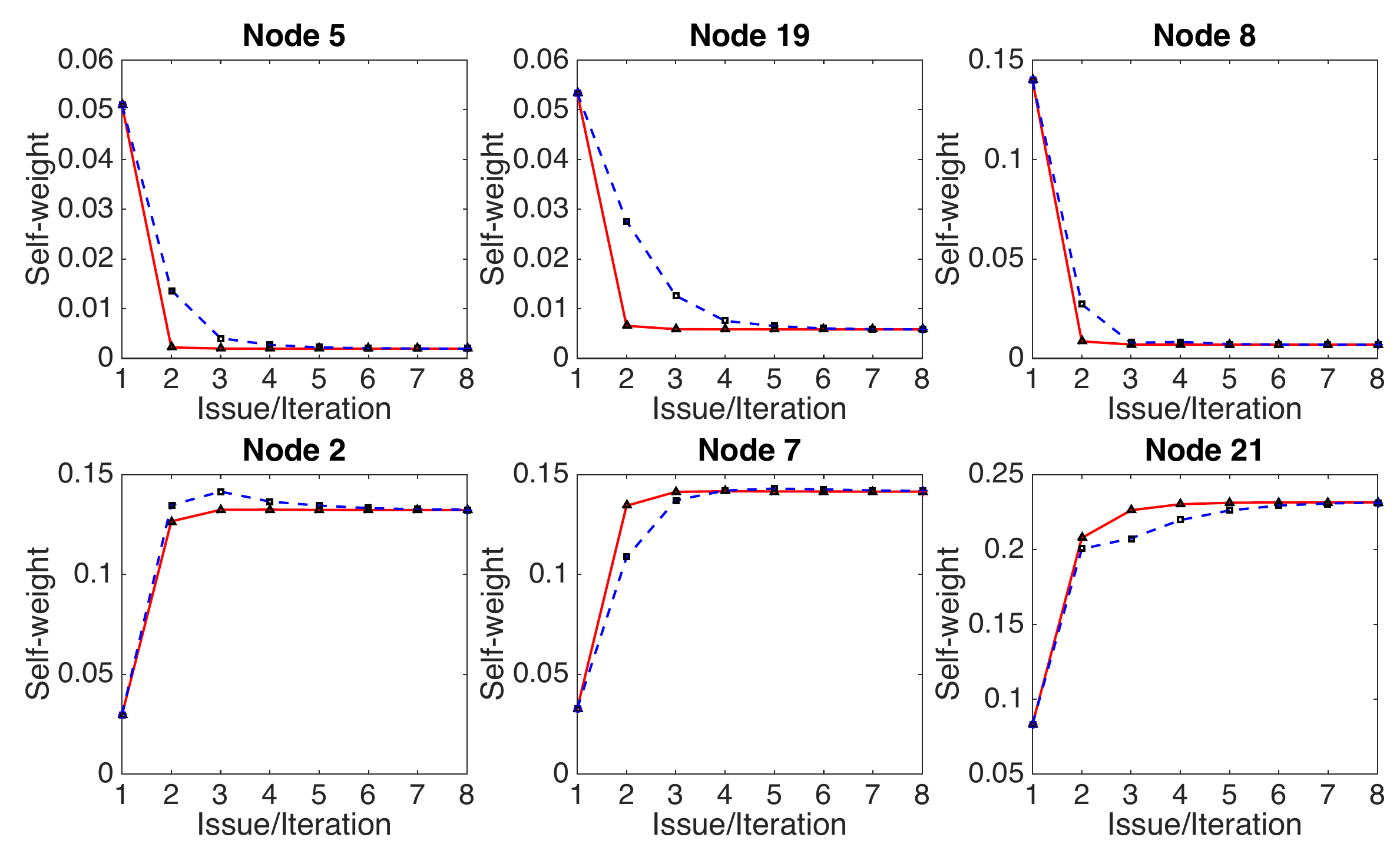}
    \caption{Self-weight evolution for the reduced Krackhardt's advice network: 
    We display the trajectories of 6 nodes with the same initial condition; the 
    dot lines represent the single-timescale DF dynamics and the solid lines 
    represent the DF dynamics.}
    \label{fig:red-KAN-t1}
    \end{center}
\end{figure}

The dynamical trajectories of the self-weights generated by the
single-timescale DF model (in dot lines) and by the DF model (in solid
lines) are illustrated in Figure~\ref{fig:red-KAN-t1}. Given non-autocratic initial
conditions, both models converge to the same equilibrium, which is
independent of initial conditions. These results are consistent with 
Theorem~\ref{thm:row-stochastic-general} and Lemma~\ref{lemma:more-features}. 
Moreover, we observe from this and all following simulations that the 
single-timescale DF model has less monotonic behaviors and takes 
more iterations to converge, compared with the DF model. 

%  \end{minipage}
%  \hspace{.4cm}
%  \begin{minipage}[t]{0.47\linewidth}

%------------------------------------------
\section{The single-timescale DF model over reducible influence networks}
\label{sec:dynamics-reducible}
%------------------------------------------

The analysis in the previous section assumes that the relative interaction
matrix $C$ is irreducible, i.e., the associated digraph is strongly connected
and each node is reachable by any other node in the network. In this section we
consider two different scenarios where the social influence network is not
strongly connected as $C$ is reducible. The part of work is comparable to the
DF model analysis over reducible networks as in~\cite{PJ-NEF-FB:14m}.

First, in Subsection~\ref{sec:dynamics-1st-class} the matrix $C$ is assumed 
to be reducible and its associated digraph has globally reachable nodes. One can
easily check that such a $C$ admits a unique dominant left eigenvector. The
analysis of the single-timescale DF model in this scenario is essentially similar 
to that for an irreducible matrix $C$. On one hand, given non-autocratic initial 
conditions, the equilibrium of the single-timescale DF model is identical to that of 
the DF model with the same $C$; on the other hand, given autocratic initial
conditions, the equilibrium of the single-timescale DF model is not necessarily 
the same as that of the DF model.

Second, in Subsection~\ref{sec:dynamics-2nd-class} the matrix $C$ is assumed to 
be reducible and its associated condensation digraph has multiple %aperiodic 
sinks. We then establish the existence and attractivity of the equilibria for 
the single-timescale DF dynamics with this most general setting. Different from 
the DF model which has a unique equilibrium, any partition of social power among 
the sinks is allowable at equilibrium of the single-timescale DF model here.

\subsection{Reducible relative interactions with globally reachable nodes}
\label{sec:dynamics-1st-class}
In this subsection we consider the single-timescale DF model in the setting of
reducible $C$ with globally reachable nodes. Recall that $C$ is reducible if and
only if $G(C)$ is not strongly connected. Without loss of generality, assume
that the globally reachable nodes are $\until{r}$, for $r\leq n$, and let
$G(C_r)$ be the subgraph induced by the globally reachable nodes.
One can show that there does not exist a row-stochastic matrix $C$ with zero
diagonal and with only one globally reachable node. However, if $r=1$, by 
assuming that node $1$ is the only globally reachable node, it is necessary that
$w_{11}=1$ and then $x(0)=\vect{e}_{1}$ as $W$ is row-stochastic by definition.
The single-timescale DF dynamics then converge to $x^*=x(0)$ even if $C$ is not
well defined. We therefore assume $r\geq2$ in the following. 

\begin{theorem}[Single-timescale DF behavior with reachable nodes]
  \label{thm:reducible-reachable} 
  For $n\geq r\geq 2$, consider a single-timescale DF dynamical system 
  $x(t+1)=F(x(t))$ as defined in~\eqref{eq:sys_IF} associated with a 
  relative interaction matrix $C\in\real^{n\times{n}}$ which is 
  row-stochastic, reducible and with zero diagonal. Let $\until{r}$ be 
  the globally reachable nodes of $G(C)$.  
  %Assume that the globally reachable subgraph $G(C_r)$ is aperiodic, 
  Then the set of equilibrium points of $F$ are 
  $\{\vect{e}_1,\dots,\vect{e}_n, x^*\}$, where $x^*\in \simplex{n}$ has 
  the following properties: 
  \begin{enumerate}
  \item[\rm{(i)}] \label{fact:reducible reachable 2} if $r=2$, then
    $x^*=\{(\alpha,1-\alpha, 0,\cdots, 0)^\top\}$ for any $\alpha\in[0,1]$, 
    and the self-weights $x(t)$ exponentially converge to $x^*$ given 
    a non-autocratic initial $x(0)$;

  \item[\rm{(ii)}] \label{fact:reducible reachable 3 star} if $r\geq 3$ 
    and $G(C_r)$ has star topology with the center node $1$, then 
    $x^*=\vect{e}_1$, and the self-weights $x(t)$ asymptotically converge 
    to $\vect{e}_1$ given any non-autocratic initial $x(0)$;

  \item[\rm{(iii)}] \label{fact:reducible reachable 3 non star} if $r\geq
    3$ and $G(C_r)$ does not have star topology, then $x^* \in
    \simplex{n} \setminus\{\vect{e}_1,\dots,\vect{e}_n\}$ satisfies: 1)
    $x^*_i>0$ for $i\in \until{r}$ and $x^*_j=0$ for $j\in\fromto{r+1}{n}$, 
    and 2) the ranking of the entries of $x^*$ is equal to the ranking of 
    the eigenvector centrality scores $c$; moreover, the self-weights $x(t)$ 
    exponentially converge to $x^*$ given any non-autocratic initial $x(0)$.
  \end{enumerate}
\end{theorem}

\begin{remark}[Comparison with the DF model]
 While the DF model and the single-timescale DF model have the same
 equilibrium set over irreducible networks, this is not true anymore
 for reducible networks with globally reachable nodes. By
 Theorem~\ref{thm:reducible-reachable}, all vertices of the simplex
 $\simplex{n}$, $\{\vect{e}_1,\dots,\vect{e}_n\}$ are the equilibrium
 points of the single-timescale DF dynamical system, whereas only the
 vertices corresponding to globally reachable nodes,
 $\{\vect{e}_1,\dots,\vect{e}_r\}$, are the equilibrium points of the DF
 model.  Nevertheless, the equilibrium point $x^*$ in the interior of
 $\simplex{n}$ for the single-timescale DF dynamics is identical to
 that associated with the DF model. In both models, $x^*$ is almost
 globally attractive.  %\oprocend
\end{remark}

\subsection*{Numerical examples on reducible networks with globally reachable nodes}
In the following, we simulate the single-timescale DF dynamics on the complete
Krackhardt's advice network (as shown in Figure~\ref{fig:KAN-topology}) and 
on a reducible network with star topology on its irreducible nodes.

\begin{figure}[htp]
    \centering
    \includegraphics[width=.6\textwidth]{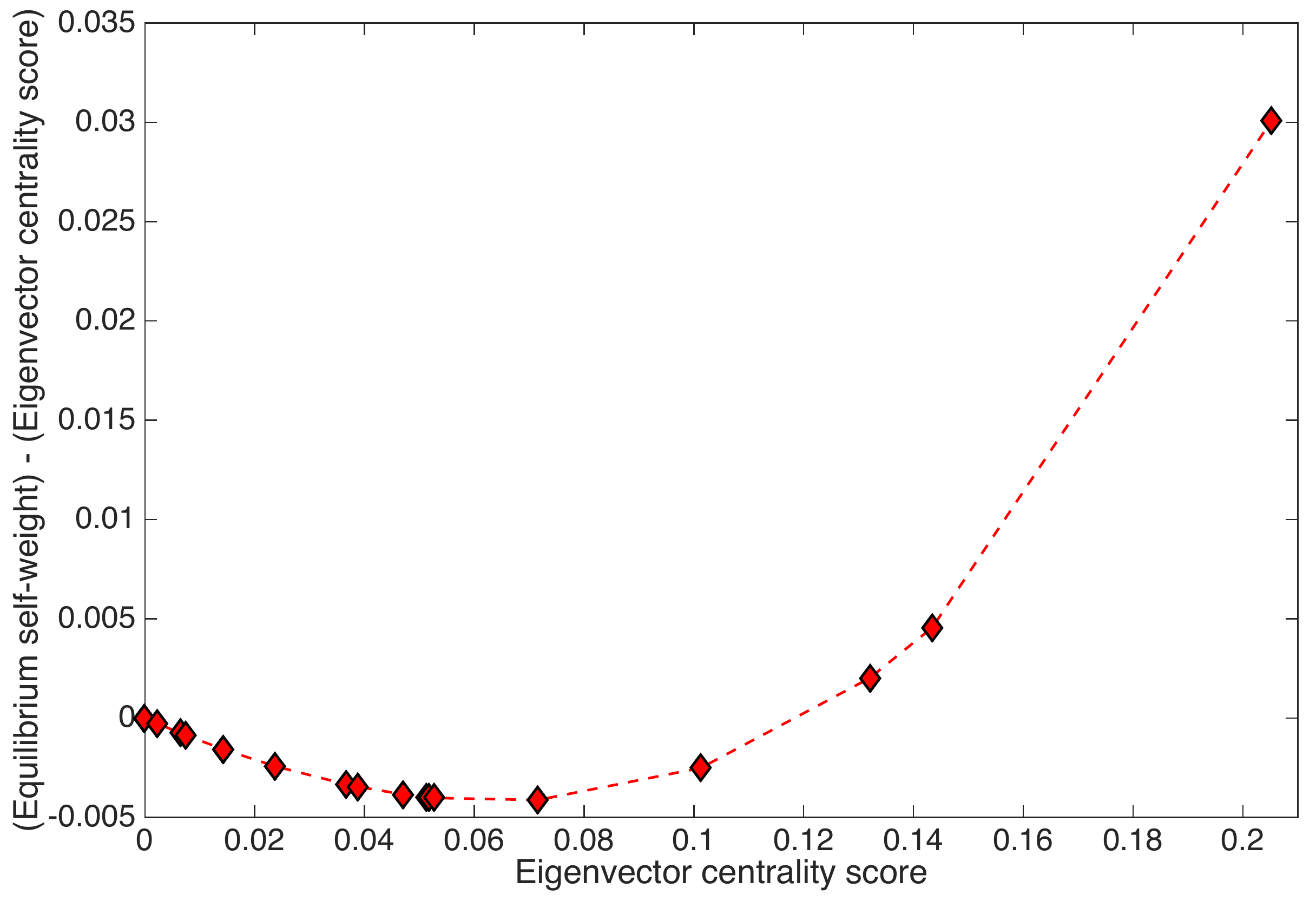}
    \caption{Comparison between the eigenvector centrality scores and
     the equilibrium self-weights for the Krackhardt's advice network: the 
     social power accumulation.}
    \label{fig:KAN-cvsx}
%  \end{minipage}
\end{figure}

\paragraph{Complete Krackhardt's advice network}
The complete Krackhardt's network, as illustrated in Figure~\ref{fig:KAN-topology}, includes four managers (i.e., individuals 
$6$, $13$, $16$ and $17$) from whom no other individual requests advice. 
Hence, this network is reducible but with globally reachable nodes (i.e., the rest $17$ individuals). Similar to 
the reduced Krackhardt's network, if individual $i$ asks for advice from $n_i$ 
different individuals, then we assume that $c_{ij}=1/n_i$ for $j$ in these $n_i$ 
individuals, and $c_{ik}=0$ for all other individuals $k$. Moreover, self-weighting 
is not considered in $C$, that is, $c_{ii}=0$ for all $i\in \until{21}$. The 
corresponding vectors $c$ and $x^*-c$ of the complete Krackhardt's advice 
network are demonstrated in Figure~\ref{fig:KAN-cvsx} to show the 
phenomenon of social power accumulation. Meanwhile, we can check that the 
ordering of the vector components of $x^*$ is consistent with that of $c$, that is, 
$x_i^*> x_j^*$ if and only if $c_i> c_j$ for $i, j \in \until{21}$.

\begin{figure}[htp]
\begin{center}
    \includegraphics[width=.8\textwidth]{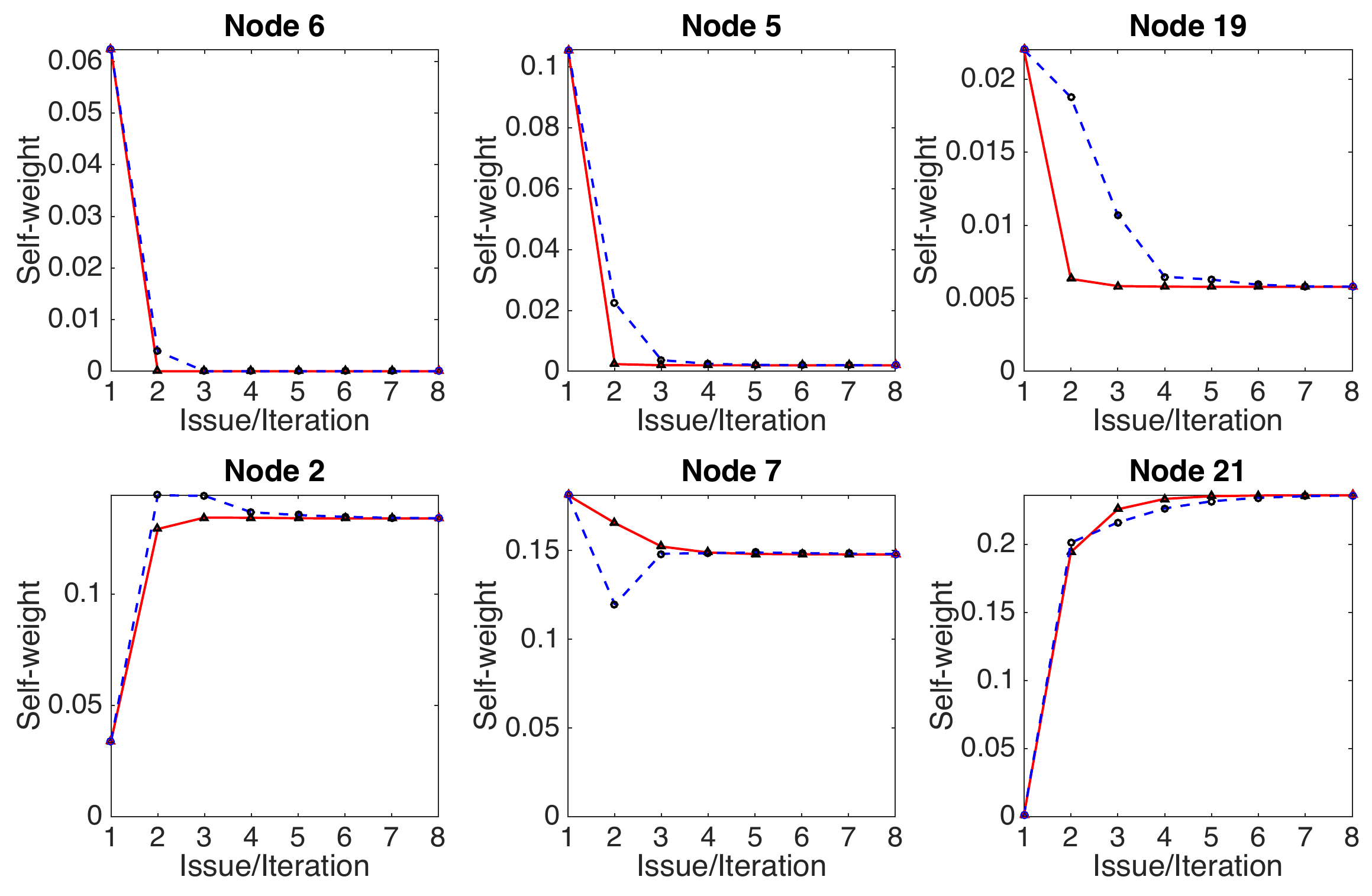}
    \caption{Self-weight evolution for the Krackhardt's advice network: we 
    display the trajectories of 6 nodes with the same initial conditions; the dot 
    lines represent the single-timescale DF dynamics and the solid lines 
    represent the DF dynamics.}
    \label{fig:KAN-t1}
    \end{center}
\end{figure}

The dynamical trajectories of the self-weights in the Krackhardt's advice network 
generated by the single-timescale DF model and the DF model are compared in Figure~\ref{fig:KAN-t1}. 
%and Figure~\ref{fig:KAN-t2}. 
For non-autocratic initial conditions, both models converge to the same 
equilibrium.

%\begin{figure}[htp]
%\begin{center}
% \vspace{-1.5cm}
%    \includegraphics[width=1\textwidth]{KN_Trj_3.pdf}
%     \vspace{-2cm}
%    \caption{Self-weight evolution for the Krackhardt's advice network: we 
%    display the trajectories of 6 nodes with different initial conditions; the dot 
%    lines represent the single-timescale DF dynamics and the solid lines 
%    represent the DF dynamics.}
%    \label{fig:KAN-t2}
%    \end{center}
%\end{figure}

\begin{figure}[htp]
\begin{center}
    \includegraphics[width=.8\textwidth]{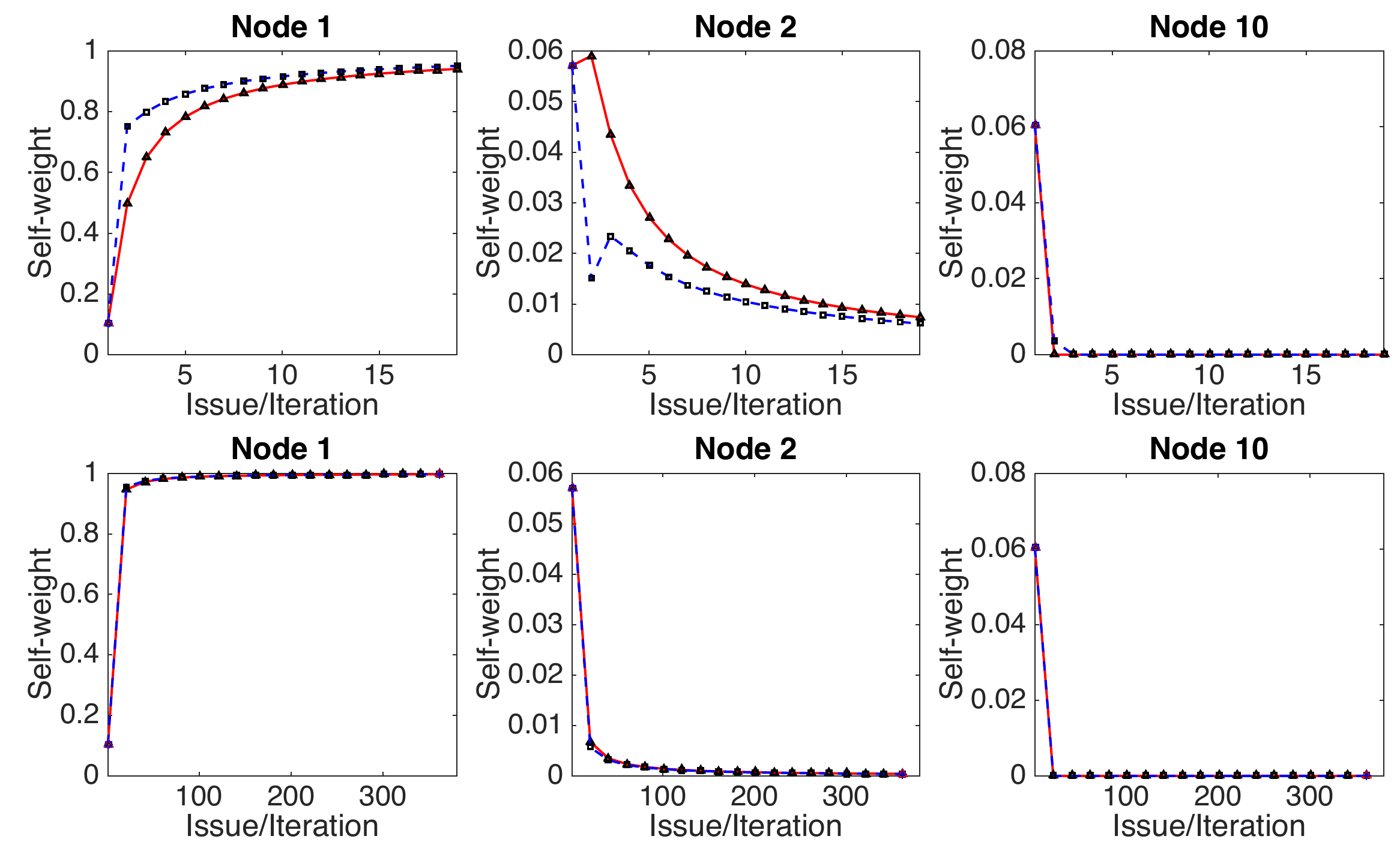}
    \caption{Self-weight evolution for a network with star topology on its 
    irreducible subgraph (that includes 
    $10$ nodes where node $10$ is reducible and node $1$ is the center): we 
    simulate both dynamics of the single-timescale DF model and of the DF model 
    with the same non-autocratic initial conditions. The dot lines represent the 
    single-timescale DF dynamics and the solid lines represent the DF 
    dynamics. The top subgraphs shows the short-term behaviors and the bottom 
    subgraphs shows the long-term behaviors. Both systems converge to the 
    same equilibrium $\vect{e}_{1}$.}
    \label{fig:rd-star-t1}
    \end{center}
\end{figure}

\begin{figure}[htp]
\begin{center}
    \includegraphics[width=.8\textwidth]{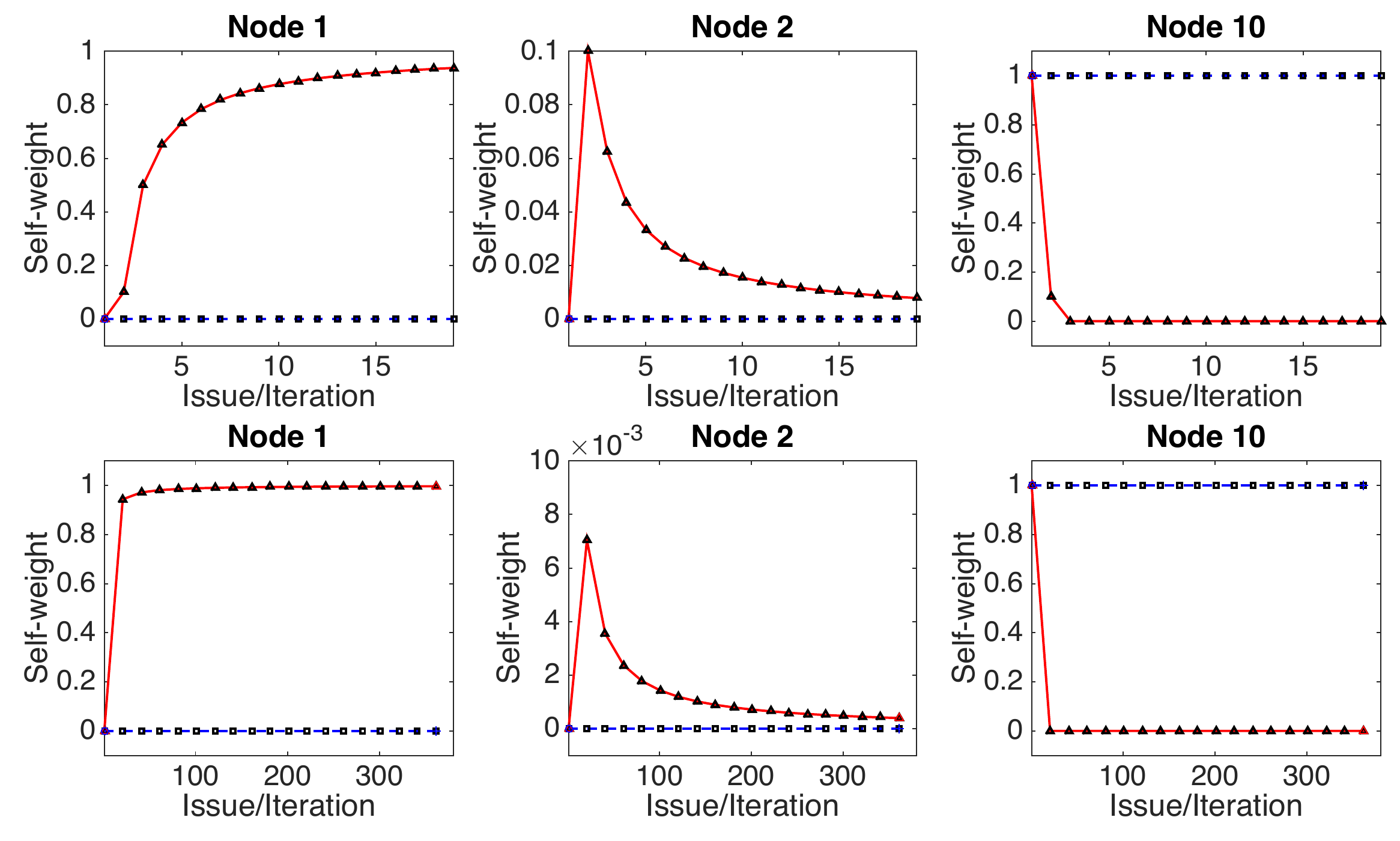}
    \caption{Self-weight evolution for the same network as in 
    Fig.~\ref{fig:rd-star-t1}: we simulate both dynamics of the single-timescale 
    DF model and of the DF model with the same autocratic initial conditions 
    $x(0)=\vect{e}_{10}$. The dot lines represent the single-timescale DF 
    dynamics and the solid lines represent the DF dynamics. The top subgraphs
    shows the short-term behaviors and the bottom subgraphs shows the 
    long-term behaviors. The two systems converge to two different equilibria 
    $\vect{e}_{10}$ and $\vect{e}_{1}$, respectively.}
    \label{fig:rd-star-t2}
    \end{center}
\end{figure}

\paragraph{A reducible network with star topology on its irreducible subgraph}
We additionally simulate the single-timescale DF dynamics on a reducible 
network with star topology on its irreducible subgraph. The single-timescale DF 
model and the DF model are compared in 
Figure~\ref{fig:rd-star-t1} and Figure~\ref{fig:rd-star-t2}. We can observe that (i) 
given a non-autocratic initial condition, both dynamical systems converge to the 
same equilibrium $\vect{e}_1$, which implies all social power is accumulated on 
individual $1$; (ii) given an autocratic initial condition on one reducible node, then 
the two systems converge to different equilibria. These statements are consistent 
with our discussion in Theorem~\ref{thm:reducible-reachable}.

\subsection{Reducible relative interactions with multiple sink components}\label{sec:dynamics-2nd-class}

In this subsection we generalize the treatment of the single-timescale DF model 
to the setting of reducible $C$ without globally reachable nodes. Such matrices 
$C$ have an associated condensation digraph $D(G(C))$ with $K\geq 2$ sinks. 
%Subject to the aperiodicity assumption on each sink, the DeGroot opinion 
%dynamical system still converges for each single issue, even though consensus is 
%not achieved for generic initial opinions.

In what follows, $n_k$ denotes the number of nodes in sink $k$,
$k\in \until{K}$, of the condensation digraph; by construction $n_k\geq 1$. 
Assume that the number of nodes in $G(C)$, not belonging to any sink 
in $D(G(C))$, is $m$, that is, $\sum_{k=1}^K n_k+m=n$. After a permutation 
of rows and columns, $C$ can be written as
\begin{equation}
C=\begin{bmatrix} C_{11}& 0 &\dots &0 &0\\
0 & C_{22} & \dots & 0&0\\
\vdots & \vdots & \ddots & \vdots &\vdots \\
0 & 0 &\dots & C_{KK} &0\\
C_{M1} & C_{M2}& \dots & C_{MK} & C_{MM}
\end{bmatrix}, \label{eq:RNR-reordered-C}
\end{equation}
where the first $(n-m)$ nodes belong to the sinks of $D(G(C))$ and the
remaining $m$ nodes do not. By construction each $C_{kk}\in
\real^{n_k\times{n_k}}$, $k\in \until{K}$, is row-stochastic and
irreducible. The Perron-Frobenius Theorem for irreducible matrices
implies that $C_{kk}$ has a unique positive dominant left eigenvector
$c_{kk}^\top=(c_{{kk}_1},\dots,c_{{kk}_{n_k}})$, satisfying
$c_{kk}\in\simplex{n_k}$, independently of whether $C_{kk}$ is
aperiodic or periodic.
Under these assumptions, the matrix $C$ has the following properties~\cite{PJ-NEF-FB:14m}: 
1) eigenvalue $1$ has geometric multiplicity equal to $K$, the number of sinks in 
the condensation digraph $D(G(C))$;
2) $C$ has $K$ dominant left eigenvectors associated with 
eigenvalue $1$, denoted by ${c^{k}}^\top \in\real^{n}$ for $k \in \until{K}$ and $c^k_i>0$ 
if and only if node $i$ belongs to sink $k$. 
We may check that $c_i^k=c_{{kk}_j}$ for $j=i-\sum_{l=1}^{k-1}n_
\ell$. We also denote $x=(x_{11}^\top, x_{22}^\top, \dots, x_{KK}^\top,x_{MM}^\top)^\top$, 
where $x_{kk}=(x_{{kk}_1},\dots,x_{{kk}_{n_k}})^\top\in \real^{n_k}$ are the 
self-weights associated with sink $k$. Similarly, $x_i=x_{{kk}_j}$ for 
$j=i-\sum_{l=1}^{k-1}n_\ell$. 
Given $x$ and $C$ with the form~\eqref{eq:RNR-reordered-C}, the corresponding
$W$ has the following form:
\begin{equation}
W=\begin{bmatrix} W_{11}& 0 &\dots &0 &0\\
0 & W_{22} & \dots & 0&0\\
\vdots & \vdots & \ddots & \vdots &\vdots \\
0 & 0 &\dots & W_{KK} &0\\
W_{M1} & W_{M2}& \dots & W_{MK} & W_{MM}
\end{bmatrix}, \label{eq:RNR-reordered-W}
\end{equation}
where $W_{Mi}=\left(I_{m}-\diag(x_{MM}) \right)C_{Mi}$ for $i<M$ and 
$W_{kk}=\diag(x_{kk})+\left(I_{n_k}-\diag(x_{kk}) \right)C_{kk}$ for $k\until{K}$.

Similar to the discussion on the single-timescale DF model with
reducible $C$ and with globally reachable nodes, for a social network
with multiple sink components and with reducible nodes, the social
power moves from the reducible nodes (by diminishing exponentially
fast) to the sinks. The social power of each sink only increases or
remains constant depending upon the initial conditions and the network
structure. The social power dynamics in each sink are similar to those
discussed in the irreducible case
Theorem~\ref{thm:row-stochastic-general}, though the total social
power of the sink is neither equal to $1$ nor constant in general.

\begin{theorem}[Single-timescale DF behavior with multiple sinks] 
\label{thm:reducible-no-reachable}
   For $n\geq 3$, consider the single-timescale DF dynamical system 
   $x(t+1)=F(x(t))$ as defined in~\eqref{eq:sys_IF} associated with 
   a relative interaction matrix $C\in\real^{n\times{n}}$. Assume that 
   the condensation digraph $D(G(C))$ contains $K\geq 2$ %aperiodic 
   sinks and that $C$ is written as in equation~\eqref{eq:RNR-reordered-C}.
   Then the following statements hold.
\begin{enumerate}

\item[\rm{(i)}]\label{fact:RNR-eq} (\textbf{Equilibrium:})
The set of equilibrium points of $F$ is the union of the set of
vertices $\{\vect{e}_1,\dots,\vect{e}_n\}$ and of the set
$\setdef{x^*=x^*(\zeta^*)\in\simplex{n}}{\zeta^*\in\simplex{K}}$,
where $\zeta^*_k$ is the total self-weight of sink $k$ and 
where $x^*$ is uniquely determined by $\zeta^*$ and has the 
following properties:
\begin{enumerate}

\item[\rm{(i.1)}]\label{fact:RNR-eq-c2-fb} if node $i$, $i\in \until{n}$, does not
  belong to any sink, then $x_i^*=0$;
		  
\item[\rm{(i.2)}]\label{fact:RNR-eq-c3-fb} if node $i$, $i\in \until{n}$, belongs
  to sink $k\in\until{K}$ and $n_k=2$, then $x^*_i=\zeta_k^*/2$ if
  $\zeta_k^*<1$, or $x^*_{kk}=(\alpha,1-\alpha)^\top$ for some
  $\alpha\in[0,1]$ if $\zeta_k^*=1$;
		
\item[\rm{(i.3)}]\label{fact:RNR-eq-c4-fb} if node $i$, $i\in \until{n}$, belongs
  to sink $k\in\until{K}$ and $n_k\geq3$, then $x^*_i>0$ if
  $\zeta_k^*>0$, or else $x^*_i=0$ if $\zeta_k^*=0$;

\item[\rm{(i.4)}]\label{fact:RNR-eq-c5-fb} for sinks with $n_k\geq3$ and
  $\zeta_k^*>0$, the ranking of the entries of the vector $x_{kk}^*$
  is equal to the ranking of the eigenvector centrality scores
  $c_{kk}$. % in the same sink $k$.
		  
\end{enumerate}

\item[\rm{(ii)}]\label{fact:RNR-constant-G} (\textbf{Monotonicity of
  sink social power:}) For all $t\geq 0$, the sink social power 
  $\zeta_{k}(t)$, equal to the sum of the  
  individual self-weights in each sink $k\in \until{K}$, is non-decreasing, 
  i.e., $\zeta_k(t+1)\geq \zeta_k(t)$; if
  $\zeta_k(0)=0$ for a sink $k$ and $x_i(0)=0$ for any reducible node
  $i$ such that there exists a direct path from $i$ to the sink $k$ in
  the associated influence network, then $\zeta_k^*=\zeta_k(t)=0$ for
  all $t\geq 0$.

\item[\rm{(iii)}]\label{fact:RNR-converge} (\textbf{Convergence of
  self-weights:}) For any initial $x(0)\in\simplex{n}\setminus
  \{\vect{e}_1,\dots,\vect{e}_n\}$, the self-weights $x(t)$
  exponentially converge to an equilibrium point $x^*$ as
  $t\to\infty$, where $x^*$ is specified as in statement~\rm{(i)}.
\end{enumerate}
\end{theorem}

\begin{remark}[Eigenvector centrality]
Similar to the DF model on reducible networks with multiple
sinks~\cite{PJ-NEF-FB:14m}, we may regard $\zeta_k^* c_{{kk}}$ as the
individual eigenvector centrality scores in sink $k$. A node has zero
eigenvector centrality score if it does not belong to any sink. When
the number of the sinks is $K\geq 2$ and $\zeta_k^*>0$ for all $k\in
\until{K}$, we have $\zeta_k^*c_{{kk}_i}<0.5$ for any sink with at
least two nodes. Consequently, the star topology in a sink does not
correspond to an equilibrium point with all sink social power on the center node of the sink, as
the eigenvector centrality score of the sink center is less than
$0.5$. Meanwhile, the social power accumulation is observed in each sink
$k$: for any individuals $i, j\in\until{n_k}$ with centrality scores
satisfying $c_{{kk}_i}>c_{{kk}_j}>0$, the social power is increasingly
accumulated in individual $i$ compared to individual $j$, that is,
$x^*_{{kk}_i}/c_{{kk}_i}>x^*_{{kk}_j}/c_{{kk}_j}$.%\oprocend
\end{remark}

\begin{remark}[Comparison with the DF model]
  For this most general case, the single-timescale DF model behaves
  very differently from the DF model: any partition of social power
  among the sinks is allowable at equilibrium of the single-timescale
  DF model, whereas the DF model has a single globally-attractive
  equilibrium, uniquely determined by $C$.  In addition, all vertices
  of the simplex $\{\vect{e}_1,\dots,\vect{e}_n\}$ are equilibrium
  points of the single-timescale DF model, but none of them is an
  equilibrium point of the DF model. % \oprocend
\end{remark}

\subsection*{Numerical examples on the Sampson's monastery network}

\begin{figure}[htp]
  \begin{center}
    \includegraphics[width=.8\linewidth]{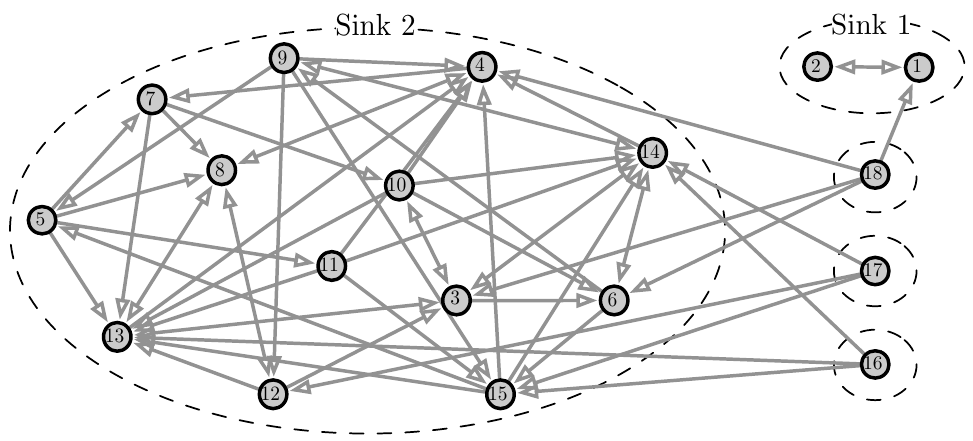}
  \end{center}
  \vspace{-.2cm}
  \caption{Sampson's monastery network}
  \label{fig:SM-topology}
\end{figure}

We demonstrate the single-timescale DF dynamics with a numerical
application to the Sampson's monastery network~\cite{SFS:69}. We
compare the single-timescale DF model with the DF model in terms 
of dynamical trajectory and equilibrium. The Sampson's monastery 
network and the corresponding $C$ have been specified in our previous 
work~\cite{PJ-NEF-FB:14m} and we use the same setup of the network. 
In particular, $C$ associated with Sampson's empirical data on esteem 
interpersonal relations is reducible. The condensation digraph 
associated with $C$ includes two sinks: sink $1$ consists of the 
nodes $\{1, 2\}$, and sink $2$ consists of the nodes $\fromto{3}{15}$, 
and the rest nodes are reducible; see Figure~\ref{fig:SM-topology}.

We simulate both the single-timescale DF model and the DF model on 
this monastery network with the same randomly selected initial states 
$x(0)\in\simplex{18}$. 
%The simulation shows that all dynamical trajectories converge to a unique 
%equilibrium self-weight vector $x^*$, given by
%\EQQ
%  x^*=  [0.0590\;\; 
%    0.0590\;\; 
%    0.1029\;\; 
%    0.2009\;\; 
%    0.0100\;\; 
%    0.0328\;\; 
%    0.0583\;\; 
%    0.1547\dots\\
%    0.0331\;\; 
%    0.0665\;\; 
%    0.0014\;\; 
%    0.0336\;\; 
%    0.1158\;\; 
%    0.0490\;\; 
%    0.0229\;\; 
%         0\;\; 
%           0\;\; 
%         0
%]^\top.
%\ENN
%Meanwhile, $\zeta_1^*=0.118$, $\zeta_2^*=0.882$, the revised eigenvector 
%centrality scores, denoted by $c^r$, can be calculated as follows:
%\EQQ
%  c^r=\zeta_1^* c^1+ \zeta_2^* c^2 =[0.0590\;\; 
%    0.0590\;\; 
%    0.1044\;\;	
%    0.1817\;\;	
%    0.0112\;\;
%    0.0359\;\;	
%    0.0622\;\;	
%    0.1479\dots\\
%    0.0363\;\;	
%    0.0702\;\;	
%    0.0016\;\;
%    0.0368\;\;	
%    0.1159\;\;	
%    0.0527\;\;	
%    0.0253\;\;	
%    0	\;\;
%    0	\;\;
%    0]^\top,
%\ENN
%and the social power accumulation threshold for sink $2$ is ${\ct}^{2}=0.1162$.

\begin{figure}[htp]
  \centering
   \includegraphics[width=.8\textwidth]{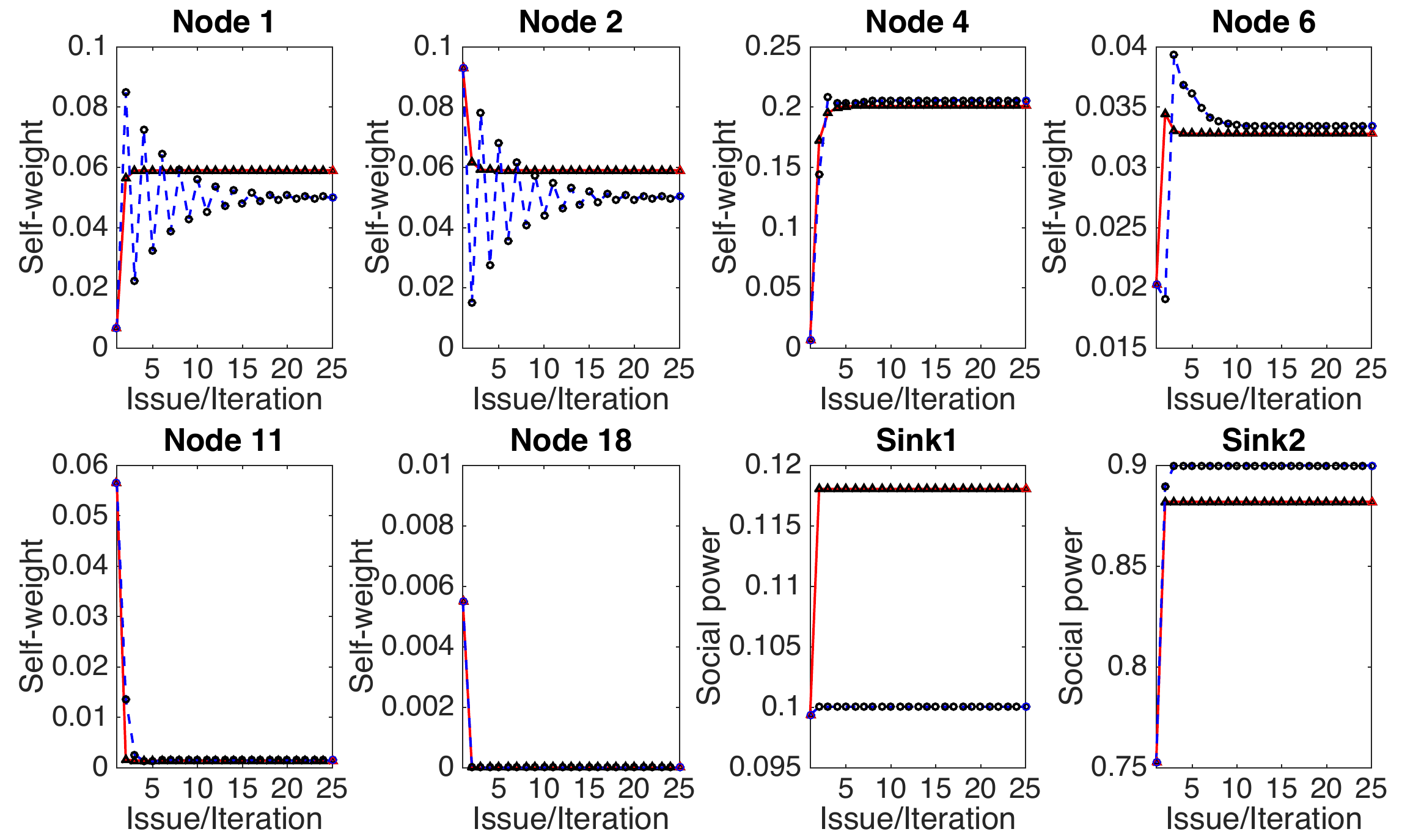}
  \caption{Self-weight evolution for the Sampson's monastery network: we
    simulate both dynamics of the single-timescale DF model and of the
    DF model with the same initial conditions. The dot
    lines represent the single-timescale DF dynamics and the solid
    lines represent the DF dynamics. We observe the
    common points for the two systems, including 1) for two nodes
    $\{1, 2\}$ in sink $1$ with $n_1=2$, the equilibrium self-weights
    are strictly positive and equal; 2) for the nodes in sink $2$ with
    $n_2=13$, all equilibrium self-weights are strictly positive and
    $x_i^*>x_j^*$ if and only $c_i^2>c_j^2$, in particular, node $4$
    has the max eigenvector centrality score in the sink, node $11$
    has the min score, and node $6$ has a score in between; 3) the
    nodes $\{16, 17, 18\}$, which do not belong to any sink, have zero
    equilibrium self-weights. We also observe the differences between
    two systems, including 1) the convergence behaviors for a sink
    with two nodes are significantly different; 2) the equilibrium
    self-weight sums are different for each sink between two systems,
    even given the same initial conditions; 3) the convergence of the
    self-weight sum at each sink occurs in two steps for the DF model, 
    but it may take more steps for the single-timescale DF model.}
  \label{fig:SM-trajectories}
 \end{figure}

The dynamical trajectories of $6$ selected nodes in the Sampson's
monastery network are
illustrated in the first $6$ subgraphs of
Figure~\ref{fig:SM-trajectories}. The trajectories of the total
self-weights in the two sinks under the same set of initial conditions are
shown in the last two subgraphs of Figure~\ref{fig:SM-trajectories}.

In addition to the differences observed from
Figure~\ref{fig:SM-trajectories}, we also note that, given different
initial conditions and a constant $C$, the DF model always converges
to the same equilibrium (see~\cite{PJ-NEF-FB:14m}), but the
single-timescale DF model converges to different equilibria by
simulation.  Specifically, regarding the DF model, the reducible nodes
have $0$ self-weights after the second issue discussion iteration and
the sum of the self-weights for each sink after the second iteration
is uniquely determined by $C$ but not $x(0)$. Moreover, the sink
social power for each sink keeps constant afterwards. Regarding the 
single-timescale DF model, the social power on reducible nodes 
converges to $0$ exponentially in general. Then at each iteration 
social power keeps migrating from reducible nodes to their connected 
sinks. Such dynamics depend not only upon $C$ but also upon the 
self-weight profile $x(t)$. As a result, each sink social power keeps 
increasing. The simulations may illustrate how different $x(0)$ 
lead to different social power evolving processes and, therefore, 
different equilibria.

%------------------------------------------
\section{Conclusion}\label{sec:conclusion}
%------------------------------------------
In this paper we have characterized the equilibrium and asymptotic
behavior of a single-timescale DF model for the evolution of social
power in a social influence network. Compared with the DF
model, a fundamental assumption in this modified model is that
individual social power evolves at the same timescale as the group
opinion forms. That is to say, social power is updated without waiting
for opinion consensus. We have derived a concise dynamical model for
the single-timescale DF evolution and completely characterized its
asymptotic properties on both irreducible and reducible networks; our
results are consistent with the partial and independent analysis
in~\cite{ZX-JL-TB:15}.  We have also compared the new model with 
the DF model in terms of their dynamical behaviors. The
analytical and numerical results show that (i) the single-timescale DF
model has the same behavior as the DF model over irreducible networks;
(ii) the single-timescale DF model behaves differently from the DF
model over reducible networks: the new model has a broader equilibrium
set including all autocratic points, and including equilibrium points
corresponding to any partition of social power among the sinks if the
underlying network has multiple sink components. Meanwhile, social
power accumulation is also observed in the new model.

This paper completes the application of reflected appraisal mechanism 
to DeGroot's opinion dynamics model and extends the validity and 
scope of the original analysis on the DF model. This paper, together 
with other efforts on time-varying influence networks, establishes 
some robustness on social power and self-appraisal evolution predicted 
by the DF model with respect to modeling uncertainties. Much work 
remains to be done in order to understand social power evolution on
various opinion formation processes. The potential examples include 
the Friedkin-Johnsen model~\cite{NEF-ECJ:99, NEF-ECJ:11}, where 
individuals tend to anchor their opinions on their initial values, and 
include influence networks with non-cooperative individuals (e.g., a 
preliminary work on existence of stubborn 
individuals~\cite{AM-PJ-NEF-FB:13x}). 

%%%------------------------------------------
%%\section*{Acknowledgments}
%%%------------------------------------------
%
%
%
%\bibliographystyle{IEEEtran}
%% \bibliographystyle{plainurl}
%% \bibliography{alias,Main,FB}

\appendices

\section{Proof of Theorem~\ref{thm:row-stochastic-star}}
\label{pf:thm:row-stochastic-star} 
\begin{IEEEproof}
Regarding fact~\rm{(i)}, we first show that the set of vertices $\{\vect{e}_1,\dots,\vect{e}_n\}$ 
are the fixed point of the dynamical system~\eqref{eq:sys_IF}. 
Given $x(t)=\vect{e}_i$ and any $C$, it is clear that 
\begin{equation*}
x(t+1)= F(\vect{e}_i)= C^\top \vect{e}_i + (I-C^\top) \vect{e}^2_i= 
C^\top \vect{e}_i + (I-C^\top) \vect{e}_i= \vect{e}_i.
\end{equation*}
Second, for $C$ associated with star topology, we show that there does not exist a fixed point in the simplex except the vertices. By contradiction, 
assume that there exists a vector $x\in\simplex{n}\setminus\{\vect{e}_1,\dots,\vect{e}_n\}$ 
such that $x=F(x)$. The fixed point equation $x=F(x)$ implies 
\begin{equation}
x_i= \sum_{j=1, j\neq i}^n c_{ji} \left( x_j - x_j^2 \right) + x_i^2, \quad \text{for all } i\in \until{n}. \label{eq:x_i}
\end{equation}
If $C$ is with star topology and the central node is $1$, then 
$c_{ij}=0$, $c_{j1}=1$ and $c_{1j}>0$ for all $j, i \in \fromto{2}{n}$. 
Especially, $c_{1j}>0$ for all $j\in \fromto{2}{n}$ because, otherwise, 
$C$ is reducible as $c_{ij}=0$ for all $i\in\until{n}$ given $j$. Therefore, 
from~\eqref{eq:x_i},
\begin{equation*}
\begin{split}
& x_j=  c_{1j} \left( x_{1} - x_{1}^2 \right) + x_j^2, \quad \text{for all }j\in \fromto{2}{n}, \\
& x_1= \sum_{j=2}^n\left( x_j - x_j^2 \right) + x_1^2.
\end{split}
\end{equation*}
That is to say, 
\begin{equation}
x_1-x_1^2= \sum_{j=2}^n\left( x_j - x_j^2 \right),\label{eq:star-e1}
\end{equation}
which implies $x_1-x_1^2>0$ as $x\in\simplex{n}\setminus\{\vect{e}_1,\dots,\vect{e}_n\}$,
and hence $(x_j-x_j^2)=c_{1j} \left( x_{1} - x_{1}^2 \right) >0$ for all $j\in \until{n}$ 
as $c_{1j}>0$. Moreover, as $x_i(1-x_i)$ is concave for $x_i\in [0, 1]$, given $n\geq 3$ and $x \in \interior{\simplex{n}}$, 
we have 
\begin{equation}
\sum_{j=2}^n x_j(1-x_j)> \sum_{j=2}^n \frac{x_j}{1-x_1} (1-x_1)x_1=x_1(1-x_1),\label{eq:star-e11}
\end{equation} 
which contradicts equation~\eqref{eq:star-e1}. Overall, for $C$ with star 
topology, all fixed points of the dynamical system~\eqref{eq:sys_IF} are 
the vertices of the simplex. 
%Note that $n\geq 3$ and $x_j>0$ for $j\in \until{n}$ together imply $\frac{x_j}{1-x_1}<1$. Also note that $f(z)=z(1-z)$ is a concave function for $z\in[0,1]$ so that $f(az)>af(z)$ for all $0<a<1$. Therefore, for $a=\frac{x_j}{1-x_1}<1$ and $z=1-x_1$, we have, for $j\in\until{n}, j\neq k$,
%  \begin{equation*} 
%    f(az)=x_j(1-x_j)> \frac{x_j}{1-x_1} (1-x_1)x_1=af(z).
%  \end{equation*}
%That is to say, $\sum_{j=1, j\neq k}^n x_j(1-x_j)> \sum_{j=1, j\neq k}^n \frac{x_j}{1-x_1} (1-x_1)x_1=x_1(1-x_1)$, which is in contradiction with equation~\eqref{eq:star-e1}.

Regarding fact~\rm{(ii)}, based upon the analysis above, for $C$ with 
star topology, the dynamical system $x(t+1)=F(x(t))$ is specified as follows:
\begin{equation}
\begin{split}
& x_j(t+1)=  c_{1j} \left( x_{1}(t) - x_{1}(t)^2 \right) + x_j(t)^2, \quad \text{for all }j\in \fromto{2}{n}, \\
& x_1(t+1)= \sum_{j=2}^n\left( x_j(t) - x_j(t)^2 \right) + x_1(t)^2. \label{eq:star_e2}
\end{split}
\end{equation}
It is clear that the function $F(x)$ is continuous for $x\in \simplex{n}$. 
If $x(0)\in\simplex{n}\setminus\{\vect{e}_1,\dots,\vect{e}_n\}$, then 
there exists a node $j$ such that $1>x_j(0)>0$, which, together with~\eqref{eq:star_e2}, 
implies $x_1(1)>0$. If $x_1(1)=1$, then $x(t)=\vect{e}_1$ for all $t\geq 1$, 
and if $x_1(1)<1$, then $x_j(2)>0$ for all $j\in \fromto{2}{n}$. Iteratively, 
we can show either $x(t)=\vect{e}_1$ or $x(t)>0$ for all $t\geq 1$. Moreover, 
if $x(t)>0$, then from~\eqref{eq:star-e11},
\begin{equation}
\label{eq:star_e3}
 x_1(t+1)-x_1(t)= \sum_{j=2}^n\left( x_j(t) - x_j(t)^2 \right) -(x_1(t)-x_1(t)^2) >0.
\end{equation}

Define a Lyapunov function candidate $V(x)=1-x_1$ for
$x\in\simplex{n}$. A sublevel set of $V$ is defined as $\left\{x \mid
V(x) \leq \beta \right\}$ for a given constant $\beta$. It is clear
that 1) any sublevel set of $V$ is compact and invariant, 2) $V$ is
strictly decreasing anywhere along the trajectory of $x(t)$ in
$\simplex{n}\setminus\{\vect{e}_1,\dots,\vect{e}_n\}$, and 3) $V$ and
$F$ are continuous. Therefore, every trajectory starting in
$\simplex{n}\setminus\{\vect{e}_1,\dots,\vect{e}_n\}$ converges
asymptotically to the equilibrium point $\vect{e}_1$ by the Lyapunov
theorem for discrete-time dynamical systems.  
\hfill \end{IEEEproof}

\section{Proof of Theorem~\ref{thm:row-stochastic-general}}
\label{pf:thm:row-stochastic-general} 
\begin{IEEEproof}
Regarding fact~\rm{(i)}, the equilibria of the influence evolution
system~\eqref{eq:sys_IF} include all vertices of the simplex as we
already demonstrate in Theorem~\ref{thm:row-stochastic-star}. Now, for
$C$ irreducible and without star topology we show that there exists a
unique $x^*\in \interior{\simplex{n}}$ satisfying $x^*=F(x^*)$ and
that the ordering of the elements of $x^*$ is consistent with that of
$c$. The fixed points of the dynamical system~\eqref{eq:sys_IF} shall
satisfy 
\begin{equation} 
\label{eq:thm-e0} 
x^*-{x^*}^2= C^\top (x^*-{x^*}^2).  
\end{equation} 
It is clear that
if $x^*\notin\{\vect{e}_1,\dots,\vect{e}_n\}$, then $x^*-{x^*}^2\neq
0$. Therefore, $(x^*-{x^*}^2)$ is a scalar multiple of the left
eigenvector of $C$ associated with eigenvalue $1$.  For $C$ without
star topology, we have 
\begin{equation*} 
 x^*-{x^*}^2= \alpha^* c, \quad \mbox{or
  equivalently,} \quad x^*_i=\alpha^* \frac{ c_i}{1-{x^*_i}}, \quad
\text{for all }i\in \until{n}, 
\end{equation*} 
where the scalar $\alpha^*$ is
such that $x^*\in \simplex{n}$, that is to say, 
\begin{equation*} 
\alpha^*=\frac{1}{\sum_{j=1}^n c_j/(1-x_j^*)}.  
\end{equation*} 
It is clear that
such an $x^*$ is exactly the same as the non-vertex fixed point we
obtained from the DF model. Therefore, the uniqueness of
$x^*$ is directly from Theorem~4.1 in~\cite{PJ-AM-NEF-FB:13d}.

Regarding fact~\rm{(ii)}, from~\eqref{eq:model-1}, we have
\begin{equation*} 
x(t+1)=\prod_{k=0}^t W(t-k)^\top x(0),
\end{equation*} 
where $W(t-k):=W(x(t-k))$ for simplicity. If we can show the product 
$\prod_{k=0}^t W(k)$ converges, then $x(t)$ also converges. To do so, 
we claim: 
\begin{enumerate}[label={(A\arabic*)},itemindent=.5cm]
 \item for any $x(0)\in\simplex{n}\setminus\{\vect{e}_1,\dots,\vect{e}_n\}$, 
$W(t)$ is aperiodic and irreducible for all $t\geq 0$ and $x(t)>0$ for all $t\geq n-1$; 
\item the minimum positive entries of $W(t)$ are lower bounded uniformly for all $t$. 
\end{enumerate}
These two claims guarantee the exponential convergence to $x^*$ for 
the dynamical system~\eqref{eq:model-1} and~\eqref{eq:sys_IF}. 
(See Lemma D.1 in~\cite{PJ-NEF-FB:13n}.)

Regarding the first claim (A1), 
as $x(0)\in\simplex{n}\setminus\{\vect{e}_1,\dots,\vect{e}_n\}$, 
there exist $m\geq 2$ nodes  with non-zero initial self-weights. 
Without loss of generality, we assume $x_i(0)>0$ for $i\in \fromto{1}{m}$ 
and the rest $n-m$ nodes with zero initial self-weights. Then, we obtain
\begin{equation} 
 x(1)=C^\top \left(x(0)-x^2(0)\right) + x^2(0)= C^\top \begin{bmatrix}x_1(0)-x_1(0)^2 \\  
 \vdots \\x_m(0)- x_m(0)^2 \\ 0 \\ \vdots \end{bmatrix}  
 + \begin{bmatrix}x_1(0)^2 \\  \vdots \\x_m(0)^2 \\ 0 \\ \vdots \end{bmatrix}. \label{eq:thm-e1}
\end{equation} 
Since $x_i(0)<1$ for $i\in \fromto{1}{m}$, $x_i(0)-x_i(0)^2>1$. 
Moreover, since $C$ is irreducible, there exist at least on edge 
from the last $n-m$ agents to the first $m$ agents, which implies 
at least one $c_{ij}>0$ for $i>m$ and $j\leq m$. Consequently, 
based upon~\eqref{eq:thm-e1}, $x_i(1)>0$ for such $i>m$ and 
$x_k(1)>0$ for all $k\leq m$. By iteration, we obtain that $x(t)>0$ 
for all $t\geq r$ given any non-vertex $x(0)$, where $r$ is the 
diameter of the digraph associated to $C$ (i.e., the maximum 
distance between any two nodes in $G(C)$).

Furthermore, consider $W(x(0)) = \diag(x(0))+(I_n-\diag(x(0)))C$. 
Since $I_n-\diag(x(0))$ has all positive diagonal entries for 
non-vertex $x(0)$, $W(x(0))$ is irreducible. As $\diag(x(0))\neq 0$,  
$W(x(0))$ is then aperiodic and primitive. The row stochasticity 
of $W(x)$ is directly from the row stochasticity assumption on $C$.  

Regarding the second claim (A2), by the definition of $W(t)$ 
in~\eqref{def:decomposition} and the constant non-negative $C$, 
the minimum positive entries of $W(t)$ are lower bounded uniformly 
if there exists a finite time $\tau\geq 0$ such that all entries 
of $x(t)$ are lower bounded uniformly for all $t\geq \tau$. 

First, we have proved above that $x(t)>0$ for all time $t\geq r$ with
$r$ as the diameter of the digraph associated to $C$.  
%Given 
%$x(0)\in\simplex{n}\setminus\{\vect{e}_1,\dots,\vect{e}_n\}$, $0\leq x_i(0)<1$ 
%for all $i\in \until{n}$. From the single-timescale DF model, 
%  \begin{equation*}
%    x(t+1)= C^\top x(t) + (I-C^\top)x^2(t)=x^2(t)+ C^\top \left(x(t)- x^2(t)\right).
%  \end{equation*}
% That means that if $x_i(t)>0$ then $x_i(t+k)>0$ for all $k\in \natural$. Moreover, since
%   \begin{equation}
%   \label{eq:IF-model-entries}
%    x_i(t+1)= x_i(t)^2+ \sum_{j=1, j\neq i}^n c_{ji} \left( x_{j}(t) - x_{j}(t)^2 \right),
%  \end{equation}
% any neighbor of individual $i$ with $x_j(t)>0$ implies $x_i(t+1)>0$. 
% These three facts, together with a irreducible $C$, indicate that 
% $x(t)>0$ for all $t\geq  r$, where $r$ is the diameter of the digraph 
% associated to $C$ (i.e., the maximum distance between any two nodes 
% in $G(C)$).  

Second, we will show that all entries of $x(t)$ are uniformly lower 
bounded away from $0$ for all $t\geq \tau$ with some $\tau\geq 0$. 
Let $\beta:=\max_{1\leq i,j\leq n}c_{ij}$ and $x_i(t)=1-\alpha$. 
It is clear that $0< \alpha <1$ and $\frac{1}{n-1}\leq \beta\leq 1$. 
Two cases (B1) $\beta<1$ and (B2) $\beta=1$ are considered in the 
following.

If (B1) $\beta<1$, as 
  \begin{equation*} 
  % \label{eq:IF-model-entries}
    x_i(t+1)= x_i(t)^2+ \sum_{j=1, j\neq i}^n c_{ji} \left( x_{j}(t) - x_{j}(t)^2 \right),
  \end{equation*} 
we have
    \begin{equation} 
     \begin{split}
    x_i(t+1)
    %&=& x_i(t)^2+ \sum_{j=1, j\neq i}^n c_{ji} \left( x_{j}(t) - x_{j}(t)^2 \right)\nnum \\
    &\leq  x_i(t)^2+ (n-1)\beta \left(\frac{\alpha}{n-1}-\frac{\alpha^2}{(n-1)^2} \right) \\
    &= (1- \alpha)^2 + \beta \alpha -\beta \frac{\alpha^2}{(n-1)},
    \label{eq:general-proof-e1}
    \end{split}
   \end{equation} 
where the inequality holds as $\beta \geq c_{ij}$ for all $1\leq i,j\leq n$ 
and the scalar function $y-y^2$ is concave on $(0, 1)$. From~\eqref{eq:general-proof-e1}, 
if $\alpha< \frac{1-\beta}{1-\beta/(n-1)}$ or equivalently 
$\beta < \frac{\alpha-\alpha^2}{\alpha-\alpha^2/(n-1)}$, by simple calculation,  
we have $x_i(t+1)< (1-\alpha)^2 + \alpha- \alpha^2 = 1-\alpha=x_i(t)$. 
That is to say, if $x_i(t)> 1-\frac{1-\beta}{1-\beta/(n-1)}$, then $x_i(t+1)<x_i(t)$.
Moreover, $x_i(t)-x_i(t+1)\geq (\alpha -\alpha^2) - \beta \alpha +\beta \frac{\alpha^2}{(n-1)}$: when $\alpha< \frac{1-\beta}{1-\beta/(n-1)}$, the right hand of this inequality has the minimum positive value at the largest $x_i(t)$ (corresponding the smallest $\alpha$) or at the point $x_i(t)=1-\frac{1-\beta}{1-\beta/(n-1)}$; in both cases $x_i(t)-x_i(t+1)$ is strictly greater than $0$. That implies the uniform decrease of $x_i(t)$ along $t$ for $\alpha< \frac{1-\beta}{1-\beta/(n-1)}$.
 
Furthermore, if $x_i(t)=1-\alpha< 1-\frac{1-\beta}{1-\beta/(n-1)}$ with 
$\alpha> \frac{1-\beta}{1-\beta/(n-1)}$, from~\eqref{eq:general-proof-e1}, 
  \begin{equation}
  \begin{split} 
    x_i(t+1)&\leq  (1- \alpha)^2 + \beta \alpha -\beta \frac{\alpha^2}{(n-1)}\\
   &= 1+ (\beta -2) \alpha +  \Big(1-\frac{\beta}{n-1}\Big)\alpha^2\\
   &=  1+ (\beta -2) \left(\frac{1-\beta}{1-\beta/(n-1)}+b\right) +  
   \Big(1-\frac{\beta}{n-1}\Big) \left(\frac{1-\beta}{1-\beta/(n-1)}+b\right)^2\\
   %&=&1-\frac{1-\beta}{1-\beta/(n-1)} + (\beta -2) b + \Big(1-\frac{\beta}{n-1}\Big) (b^2+2b \frac{1-\beta}{1-\beta/(n-1)}), 
   &=1-\frac{1-\beta}{1-\beta/(n-1)} -\beta b+ \Big(1-\frac{\beta}{n-1}\Big) b^2,
   \label{eq:general-proof-e2}
   \end{split}
    \end{equation} 
where $b=\alpha-\frac{1-\beta}{1-\beta/(n-1)}$. It is clear that $0<b< \frac{\beta}{1-\beta/(n-1)}$. 
%Hence, the part of the right hand side of~\eqref{eq:general-proof-e2} satisfies
%\EQQ
%(\beta -2) b + \Big(1-\frac{\beta}{n-1}\Big) (b^2+2b \frac{1-\beta}{1-\beta/(n-1)})
%&=& (\beta -2) b + \Big(1-\frac{\beta}{n-1}\Big) b^2+  2(1-\beta)b\nnum\\
%&=& -\beta b+ \Big(1-\frac{\beta}{n-1}\Big) b^2\nnum\\
%&=&b \left( -\beta + \Big(1-\frac{\beta}{n-1}\Big) b\right)<0.
%\ENN
%That is to say, $x_i(t+1)<1-\frac{1-\beta}{1-\beta/(n-1)}$
Consequently, the part of the right hand side of~\eqref{eq:general-proof-e2} satisfies
 \begin{equation} 
\label{eq:general-proof-e3}
 -\beta b+ \Big(1-\frac{\beta}{n-1}\Big) b^2
=b \left( -\beta + \Big(1-\frac{\beta}{n-1}\Big) b\right)<0.
\end{equation} 
Hence, $x_i(t+1)<1-\frac{1-\beta}{1-\beta/(n-1)}$ from~\eqref{eq:general-proof-e2} and~\eqref{eq:general-proof-e3}. 
Overall, if one entry of $x(t)$ is greater than $1-\frac{1-\beta}{1-\beta/(n-1)}$, then via the single-timescale DF 
model~\eqref{eq:sys_IF}, the value of the underlying entry is uniformly decreasing until it is less than 
$1-\frac{1-\beta}{1-\beta/(n-1)}$. If one entry of $x(t)$ is less than $1-\frac{1-\beta}{1-\beta/(n-1)}$, then it 
is less than $1-\frac{1-\beta}{1-\beta/(n-1)}$ for all following iterations $(t+k)$, $k\in \natural$. In other words, 
there exists a finite time $\tau$ such that all entries of $x(t)$ for all $t\geq \tau$ are bounded away from $1$ 
uniformly. Consequently, from the equation~\eqref{eq:sys_IF}, the facts $x(t)>0$ and $C$ irreducible, we have 
all entries of $x(t)$ are also bounded away from $0$.

If (B2) $\beta=1$, without loss of generality, assume $c_{i1}=\beta=1$ for some $2\leq i\leq n$. This implies that 
the $i$-th individual only accords relative interpersonal weight to the first individual in the group. As $C$
is row-stochastic, $c_{ij}=0$ for all $2\leq j\leq n$. Moreover, as $C$ is not with star topology, at least one 
individual $j$ has $c_{j1}<1$.In the following, we will show that, for a sufficiently large $\alpha$ satisfying 
$0<\alpha<1$, if $x_1(t)>\alpha$ for $t\geq n-1$, then $x_1(t+1)\leq x_1(t)$. 

Here we first consider two exclusive and complete scenarios for the case (B2): 
\begin{enumerate}[label={(C\arabic*)},itemindent=.5cm]
\item $C$ satisfies $c_{n1}<1$ and $c_{i1}=1$ for all rest individuals $i\neq 1$; and  
\item $C$ satisfies $c_{j1}<1$ for $n\geq j>m$ and $c_{i1}=1$ for $m\geq i>1$ where $n-1>m>1$. 
\end{enumerate}
Note that 1) in scenario (C1), $c_{in}=0$ for all $i\neq 1$ and $0<c_{1n}<1$;  2) we can always 
re-arrange the indices of individuals such that scenario (C2) occurs for more that one individuals 
only accord interpersonal weights to the first individual. 

Regarding the scenario (C1), by~\eqref{eq:sys_IF} and by the fact that $C$ is not with star topology,
\begin{equation} 
\label{eq:general-proof-e4}
\begin{split}
x_1(t+1)&=\sum_{i=2}^n c_{i1} (x_i(t)-x_i(t)^2) +x_1(t)^2 \\ 
&= \sum_{i=2}^{n-1}  (x_i(t)-x_i(t)^2) + c_{n1} (x_n(t)-x_n(t)^2)+x_1(t)^2.
\end{split}
\end{equation} 
Here we also assume $c_{n1}<1$ without loss of generality. 

We have proved that $x(t)>0$ for $t\geq  r$, where $r$ is the diameter of the digraph associated to 
$C$ and $r\leq n-1$. Following the equation~\eqref{eq:general-proof-e4}, to 
prove $x_1(t+1)\leq x_1(t)$ for $x_1(t)$ sufficiently close to $1$, it is sufficient to show that 
\begin{equation} 
\label{eq:general-proof-e5}
\sum_{i=2}^{n-1}  (x_i(t)-x_i(t)^2) + c_{n1} (x_n(t)-x_n(t)^2)+x_1(t)^2
\leq x_1(t).
\end{equation}
From the dynamical system~\eqref{eq:sys_IF}, we have 
\begin{equation*}
x(t+1)-x(t)=(C^\top-I) (x(t)-x(t)^2) \quad \text{for all } t\geq 0,
\end{equation*}
or equivalently for $x(t)>0$ and $x(t)\neq x^*$,
\begin{equation*}
x(t+1) \diag(x(t))^{-1}=(C^\top-I) (\vectorones[n]-x(t))+\vectorones[n]= C^\top (\vectorones[n]-x(t))+ x(t).
\end{equation*}
Therefore, for $x_n(t)>0$ and $x_n(t)\neq x_n^*$, 
\begin{equation}
\label{eq:general-proof-e5.1}
\frac{x_n(t+1)}{x_n(t)}=\sum_{i=1}^{n-1} c_{in} (1-x_i(t)) +x_n(t)= c_{1n} (1-x_1(t)) +x_n(t)<1.
\end{equation}
That is to say, $x_n(t+1)<x_n(t)$. Moreover, as $x_n(t+1)> c_{1n} (x_1(t)-x_1(t)^2)$, the following statement also holds: 
\begin{equation} 
x_n(t)> c_{1n}(x_1(t)-x_1(t)^2)\geq \gamma (1-x_1(t))=\gamma \sum_{i=2}^n x_i(t),
\label{eq:general-proof-e6} 
\end{equation}
with $0<\gamma \leq c_{1n} x_1(t)<1$.

Moreover, based upon~\eqref{eq:general-proof-e6} and for a sufficient large $x_1(t)<1$, 
we have the following statements related:
\begin{equation}
\label{eq:general-proof-e7}
\begin{split}
&\sum_{i=2}^{n-1}  (x_i(t)-x_i(t)^2) + c_{n1} (x_n(t)-x_n(t)^2)+x_1(t)^2
< x_1(t)\\
\iff& \sum_{i=2}^{n-1}  (x_i(t)-x_i(t)^2) + c_{n1} (x_n(t)-x_n(t)^2) < \sum_{i=2}^{n}  x_i(t) - (\sum_{i=2}^{n}  x_i(t))^2 \\
\iff& (1-c_{n1}) (x_n(t)-x_n(t)^2) > (\sum_{i=2}^{n}  x_i(t))^2 - \sum_{i=2}^{n}  x_i(t)^2\\
\Longleftarrow\;\;& (1-c_{n1}) (\gamma \sum_{i=2}^n x_i(t) -\gamma^2 (\sum_{i=2}^n x_i(t))^2)\geq (\sum_{i=2}^{n}  x_i(t))^2 - \sum_{i=2}^{n}  x_i(t)^2\\
\iff& (1-c_{n1}) \gamma \sum_{i=2}^n x_i(t) \geq ((1-c_{n1}) \gamma^2+1)(\sum_{i=2}^{n}  x_i(t))^2- \sum_{i=2}^{n}  x_i(t)^2\\
\Longleftarrow\;\;& (1-c_{n1}) \gamma \sum_{i=2}^n x_i(t) \geq ((1-c_{n1}) \gamma^2+\frac{n-1}{n})(\sum_{i=2}^{n}  x_i(t))^2\\
\iff& \frac{(1-c_{n1}) \gamma}{(1-c_{n1}) \gamma^2+\frac{n-1}{n}}\geq  \sum_{i=2}^{n}  x_i(t)=1-x_1(t).
\end{split}
\end{equation}
The last statement holds for $x_1(t)\geq 1-\frac{(1-c_{n1}) \gamma}{(1-c_{n1}) \gamma^2+\frac{n-1}{n}}$, 
where $\gamma<(n-1)/n$  guarantees $0<\frac{(1-c_{n1}) \gamma}{(1-c_{n1}) \gamma^2+\frac{n-1}{n}}<1$. 
Therefore, the inequality~\eqref{eq:general-proof-e5} holds. That is, 
$x_1(t+1)\leq x_1(t)$ for $x_1(t)\geq 1-\frac{(1-c_{n1}) \gamma}{(1-c_{n1}) \gamma^2+\frac{n-1}{n}}$.
In addition, for the system~\eqref{eq:sys_IF}, we have for all $x_1(t)\leq 1-\frac{(1-c_{n1}) \gamma}{(1-c_{n1}) \gamma^2+\frac{n-1}{n}}:=\beta$,
\begin{equation} 
\label{eq:general-proof-e7.1}
x_1(t+1)= \sum_{i=2}^{n} c_{i1} (x_i(t)-x_i(t)^2) +x_1(t)^2 < \beta^2+ (1-\beta) (1-\frac{1-\beta}{n-1}) <1.
\end{equation}
This implies that there exists a finite time $\tau$ such that $x_1(t)$ for all $t\geq \tau$ are bounded away from $1$ and bounded 
way from $0$. 
%By using similar arguments as in the case (B1), it is straightforward to check the inequality $x_i(t)\geq x_i(t+1)$ holds for all the 
%rest individuals $2\leq i\leq n$ if $x_i(t)$ is sufficiently large, and \eqref{eq:general-proof-e7.1} always holds, given $C$ defined 
%in the scenario (C1). Hence, there exists a finite time $\tau$ such that all entries of $x(t)$ for all $t\geq \tau$ are bounded away 
%from $1$ and from $0$. %We may check that $x_1(t+1)=x_1(t)$ only if $x(t)=x^*$.

Regarding the scenario (C2), we may regard the set of individuals $(m+1, m+2, \cdots, n)$ as a single "node", 
as they are only directly connected to the first individual but not the rest set of individuals $(2, 3,\cdots, m)$. 
Similar arguments as for the scenario (C1) hold here to prove $x_1(t+1)\leq x_1(t)$. First, we have the similar 
statement to \eqref{eq:general-proof-e5.1}.
For any $x_j(t)>0, m\leq j\leq n$, and $x_j(t)\neq x_n^*$, 
\begin{equation*}
\begin{split}
\frac{x_j(t+1)}{x_j(t)}=&\sum_{i=1}^{n-1} c_{ij} (1-x_i(t)) +x_n(t)\\
=& \sum_{i\neq j,i=m+1}^{n} c_{ij} (1-x_i(t)) +c_{1n} (1-x_1(t)) +x_n(t)<n-m-1+1=n-m.
\end{split}
\end{equation*}
That is to say, $x_j(t+1)<x_j(t) (n-m)$ for all $m+1\leq j\leq
n$. Moreover, as the digraph associated with $C$ is irreducible, there
exists at least one $m+1\leq j\leq n$ such that $x_j(t+1)> c_{1j}
(x_1(t)-x_1(t)^2)$, this implies $x_j(t)\geq \gamma_j (1-x_1(t))$ for
some $\gamma_j<1$ and independent of time $t$. Consequently, there
exists at least one different individual $m+1\leq i\leq n, i\neq j$
such that $x_i(t+1)> c_{ji}
(x_j(t)-x_j(t)^2)>c_{ji}x_j(t)>\gamma_i(1-x_1(t))$. Similarly, we have
all individuals $m+1\leq j\leq n$ with $c_{j1}<1$ satisfying
$x_j(t)\geq \gamma_j (1-x_1(t))$ for some $\gamma_j<1$.  Second, we
have the similar statement to \eqref{eq:general-proof-e7}:
\begin{equation}
\label{eq:general-proof-e8}
\begin{split}
&\sum_{i=2}^{m}  (x_i(t)-x_i(t)^2) + \sum_{j=m+1}^{n} c_{j1} (x_j(t)-x_j(t)^2)+x_1(t)^2
\leq x_1(t)\\
\iff& \sum_{i=2}^{m}  (x_i(t)-x_i(t)^2) + \sum_{j=m+1}^{n} c_{j1} (x_j(t)-x_j(t)^2) \leq  \sum_{i=2}^{n}  x_i(t) - (\sum_{i=2}^{n}  x_i(t))^2  \\
\iff& \sum_{j=m+1}^{n} (1-c_{j1}) (x_j(t)-x_j(t)^2) \geq (\sum_{i=2}^{n}  x_i(t))^2 - \sum_{i=2}^{n}  x_i(t)^2 \\
\Longleftarrow\;\; & (n-m)(1-c_{j1}) (\gamma_j \sum_{i=2}^n x_i(t) -\gamma_j^2 (\sum_{i=2}^n x_i(t))^2)\geq (\sum_{i=2}^{n}  x_i(t))^2 - \sum_{i=2}^{n}  x_i(t)^2   \\
&\quad (\mbox{where $j=\argmin_{m<i\leq n}(1-c_{i1}) (x_i(t)-x_i(t)^2)$}) \\ 
\Longleftarrow\;\; & (1-c_{j1}) \gamma_j \sum_{i=2}^n x_i(t) \geq ((1-c_{j1}) \gamma_j^2+1)(\sum_{i=2}^{n}  x_i(t))^2- \sum_{i=2}^{n}  x_i(t)^2 \\
\Longleftarrow\;\; & (1-c_{j1}) \gamma_j \sum_{i=2}^n x_i(t) \geq ((1-c_{j1}) \gamma_j^2+\frac{n-1}{n})(\sum_{i=2}^{n}  x_i(t))^2 \\
\iff&\frac{(1-c_{j1}) \gamma_j}{(1-c_{j1}) \gamma_j^2+\frac{n-1}{n}}\geq  \sum_{i=2}^{n}  x_i(t)=1-x_1(t).
\end{split}
\end{equation}
Hence, for $x_1(t)\geq 1-\frac{(1-c_{j1}) \gamma_j}{(1-c_{j1}) \gamma_j^2+\frac{n-1}{n}}$, from \eqref{eq:general-proof-e8}, we have 
\begin{equation*}
x_1(t+1)=\sum_{i=2}^{n} c_{i1} (x_i(t)-x_i(t)^2) +x_1(t)^2
\leq x_1(t).
\end{equation*}
As \eqref{eq:general-proof-e7} always holds, we can prove that there exists a finite time $\tau_1$ such that 
$x_1(t)$ for all $t\geq \tau_1$ are bounded away from $1$ and from $0$. 
%The rest proof that $x_i(t), i\neq 1,$ is always bounded way from $1$ and $0$ for all $t\geq \tau_i$ and for some finite $\tau_i$, 
%is straightforward. On the whole, given $\tau=\max(\tau_i)_{1\leq i\leq n}$, all entries of $x(t)$ are bounded away from $1$ and 
%from $0$ for all $t\geq \tau$, given $C$ defined in the scenario (C2).

Overall, given any $C$ irreducible and row-stochastic, each individual $i$ in the network must satisfy one among 
the three cases (B1) (although we assume all non-zero $c_{ij}<1$ in (B1), we only require $c_{ji}<1$ for all $j\neq i$ 
given $i$ in the proof), (C1) and (C2). That is, there always exists a finite time $\tau_1$ such that $x_1(t)$ for all 
$t\geq \tau_1$ are bounded away from $1$ and from $0$, given non-vertex $x(0)$. Consequently, from the 
equation~\eqref{eq:sys_IF} and the facts $C$ irreducible, there always exists a finite time $\tau_i$ such that all 
entries of $x_i(t)$ for all $t\geq \tau_i$ are bounded away from $1$ and from $0$ uniformly. Hence, there exists 
a finite time $\tau$ such that all entries of $x(t)$ for all $t\geq \tau$ are bounded away from $1$ and from $0$.
As a result, the claim (A2) holds, which completes the proof of fact~\rm{(ii)}.
\hfill\end{IEEEproof}

\section{Proof of Theorem~\ref{thm:reducible-reachable}}
\label{pf:thm:reducible-reachable} 
\begin{IEEEproof}
By definition, 
\begin{equation*}
x_i(t+1)= x_i(t)^2+ \sum_{j=1, j\neq i}^n c_{ji} \left( x_{j}(t) - x_{j}(t)^2 \right).
\end{equation*}
As $x(0)$ is in a simplex, if $x_i(0)=1$ then $x_j(0)=0$ for all $j\neq i$. It is 
clear that $x_i^*=x_i(1)=1$ and therefore, $x^*= \vect{e}_i$ given 
$x(0)= \vect{e}_i$ for all $i\in \until{n}$. That is to say, $\{\vect{e}_1,\dots,\vect{e}_n\}$ 
are always the fixed points of the dynamical system~\eqref{eq:sys_IF}. 

Regarding fact~\rm{(i)}, without loss of generality, we assume that node $1$ and node $2$ 
are globally reachable. The corresponding $C$ has the following block matrix form
\begin{equation}
C=\begin{bmatrix} C_1  & 0  \\
C_{21}  &  C_{22}   
\end{bmatrix}, \label{eq:condensation-graph-2g}
\end{equation}
where $C_1\in \real^{2\times 2}$ is row stochastic, and $C_{22}\in \real^{n-2\times n-2}$ 
is substochastic as $C_{21}\geq 0$ and $C_{21}\neq 0$. 
Given $x(t)\in \simplex{n} \setminus\{\vect{e}_1,\dots,\vect{e}_n\}$, the weight matrix 
$W(t)$ has the block matrix form via \eqref{def:decomposition} as follows.
\begin{equation}
W(t)=\begin{bmatrix} W_1(t)  & 0  \\
W_{21}(t)  &  W_{22}(t)   
\end{bmatrix}. \label{eq:w-2g}
\end{equation}
Here
\begin{equation*}
\begin{split}
&W_1(t):=W_1(x_{(1,2)} (t))=\diag{x_{(1,2)} (t)} + (I_2-\diag(x_{(1,2)}(t)))C_1,\\
&W_{21}(t):=W_{21}(x_{(3,\cdots,n)} (t))= (I_{n-2}-\diag(x_{(3,\cdots,n)}(t)))C_{21},\\
&W_{22}(t):=W_{22}(x_{(3,\cdots,n)} (t))=\diag{x_{(3,\cdots,n)} (t)} + (I_{n-2}-\diag(x_{(3,\cdots,n)}(t)))C_{22},
\end{split}
\end{equation*}
given $x_{(1,2)}:=\begin{bmatrix} x_1&  x_2 \end{bmatrix}^\top$ and 
$x_{(3,\cdots,n)}:=\begin{bmatrix} x_3&\cdots&  x_n \end{bmatrix}^\top$ .

The single-timescale DF dynamics associated with $C$ in~\eqref{eq:condensation-graph-2g} is as follows.
\begin{equation*}
\begin{split}
  &x_{(1,2)}(t+1)= W_{1}(t)^\top x_{(1,2)}(t)+ W_{21}(t)^\top x_{(3,\cdots,n)}(t), \\
  &x_{(3,\cdots,n)}(t+1)= W_{22}(t)^\top x_{(3,\cdots,n)}(t), \hspace{3cm} t=0, 1, 2,\dots.
  \end{split}
\end{equation*}
As $C_{22}$ is substochastic and $x(0)\in \simplex{n}\setminus\{\vect{e}_1,\dots,\vect{e}_n\}$, 
$W_{22}(0)$ is substochastic. That is, $\sum_{i=3}^n x_i(1)\leq \sum_{i=3}^n x_i(0)$ and  
$\max_{3\leq i\leq n} x_i(1)\leq \max_{3\leq i\leq n} x_i(0)$ for $\max_{3\leq i\leq n} x_i(0)\neq 0$. 
These statements hold for all $t\geq 0$ iteratively. In particular, for $x_{(3,\cdots,n)}(0)= 0$,
$x_{(3,\cdots,n)}(t)= 0$ for all $t\geq 0$. Moreover, as $C$ is reducible and has globally reachable nodes,
given any initial conditions of $x_{(3,\cdots,n)}(0)$, the zero and non-zero pattern of $x_{(3,\cdots,n)}(t)$
shall keep constant for all $t\geq n-2$. That is, for $3\leq i\leq n$ and $t\geq n-2$, if $x_{i}(t)>0$ then 
$x_{i}(t+k)>0$ for all finite $k\geq 0$, and if $x_{i}(t)=0$ then $x_{i}(t+k)=0$ for all $k\geq 0$.

Next, we will show that, given $x_{(3,\cdots,n)}(0)\neq 0$, $\lim_{t\rightarrow \infty} x_{(3,\cdots,n)}(t)=0$
exponentially. By appropriately re-indexing all individual $3\leq i\leq n$, we have $C_{22}$ have the following
normal form:
\begin{equation*}
C_{22} = \begin{bmatrix} A_{11}  &0        &0       &\cdots &0     \\ 
												 A_{21}  &A_{22}   &0       &\cdots &0     \\     
												 \vdots  &\vdots   &\vdots  &       &\vdots\\
												 A_{m1}  &A_{m2}   &A_{m2}  &\cdots &A_{mm}
					 \end{bmatrix}.
\end{equation*}
If $C_{22}$ is irreducible then $m=1$; otherwise, each block matrix $A_{ii}$ is irreducible for $i=\until{m}$. 
Moreover, as $C_{22}$ is substochastic and $C$ is stochastic and have globally reachable nodes, each $A_{ii}$
is substochastic with at least one row sum strictly less than $1$. Consequently, from~\eqref{eq:w-2g}, we have
\begin{equation*}
W_{22}(t)= \begin{bmatrix} B_{11}(t) &0         &0        &\cdots &0\\ 
													B_{21}(t) &B_{22}(t)  &0        &\cdots &0\\     
													\vdots    &\vdots    &\vdots   &       &\vdots\\
													B_{m1}(t) &B_{m2}(t) &B_{m2}(t)&\cdots &B_{mm}(t)
					 \end{bmatrix}.
\end{equation*} 
where $B_{ii}(t)=\diag(x_{s_i} (t)) + \left(I_{|s_i|}-\diag(x_{s_i}(t)) \right)A_{ii}$, 
given $s_i$ is the set of individuals corresponding to the rows evolving in the 
block matrix $A_{ii}$ and $|s_i|$ denotes the cardinality of the set $s_i$. It is clear 
that $B_{ii}(t)$ is irreducible, substochastic, and has at least one row sum 
strictly less than $1$, for all $t\geq 0$. Moreover, as the maximum of the elements 
of $x_{s_i} (t)$ is less than or equal to the maximum of the elements of 
$x_{(3,\cdots,n)} (0)$, the elements of $B_{ii}(t)$ are upper bounded uniformly for 
all $t\geq 0$. Meantime, all $B_{ii}(t)$ for $t\geq n-2$ shall have the same zero 
and non-zero pattern on elements. As a result of all these facts and from~\cite[Corollary~4.11]{FB:17}, $\prod_{k=0}^t B_{ii}(t)$ converges to $\vectorzeros[]$ 
exponentially for each block matrix and hence, $\prod_{k=0}^t W_{22}(t)$ converges to 
$\vectorzeros[]$ exponentially. From~\eqref{eq:sys_IF}, $x_{(3,\cdots,n)} (t)$ converges to $\vectorzeros[n-2]$ 
exponentially. 

As $x_{(1,2)}(t+1)= W_{1}(t)^\top x_{(1,2)}(t)+ W_{21}(t)^\top x_{(3,\cdots,n)}(t)$ and $C_1=\begin{bmatrix} 0 & 1\\ 1 & 0 \end{bmatrix}$, 
we have  
\begin{equation*}
W_1(t)=\diag{x_{(1,2)} (t)} + (I_2-\diag(x_{(1,2)}(t)))C_1=\begin{bmatrix} x_1(t) & 1-x_1(t)\\ 1-x_2(t) & x_2(t) \end{bmatrix}.
\end{equation*}
Once $x_{(3,\cdots,n)} (t)$ converges to $\vectorzeros[n-2]$ exponentially, $x_{(1,2)}(t)$ 
simultaneously converges to an equilibrium $x^*_{(1,2)}$  satisfying
\begin{equation*}
x^*_{(1,2)}:=\lim_{t\rightarrow \infty}x_{(1,2)}(t) = \lim_{t\rightarrow \infty} W_{1}(t)^\top x_{(1,2)}(t).
\end{equation*}
That is 
\begin{equation*}
\begin{bmatrix} x_1^* \\  x_2^* \end{bmatrix}=\begin{bmatrix} x_1^* & 1-x_2^*\\ 1-x_1^* & x_2^*  \end{bmatrix} \begin{bmatrix} x_1^* \\  x_2^* \end{bmatrix} =\begin{bmatrix} (x_1^*)^2 +x_2^*- (x_2^*)^2\\ x_1^*-(x_1^*)^2 +(x_2^*)^2 \end{bmatrix} \iff x_2^*- (x_2^*)^2=x_1^*- (x_1^*)^2.
\end{equation*}
As $\lim_{t\rightarrow \infty} (x_1(t)+x_2(t))=1$, $x_2^*- (x_2^*)^2=x_1^*- (x_1^*)^2$ 
holds for any pair $(x_1^*, x_2^*)$ satisfying $x_1^*+x_2^*=1$.

Regarding fact~\rm{(ii)}, the similar arguments in~\rm{(i)} can prove all $\{x_i(t)\}$ 
corresponding to reducible nodes converge to $0$ exponentially. Consequently, if 
$\sum_{i=r+1}^n x_{(i)} (t)= \beta(t)$ sufficiently small, the following statement 
similar to~\eqref{eq:star_e3} holds
\begin{equation}
\label{eq:star_r_e3}
 x_1(t+1)-x_1(t)= \sum_{j=2}^r\left( x_j(t) - x_j(t)^2 \right) -(x_1(t)-x_1(t)^2) >0.
\end{equation}
for all $x_1(t) \leq 1-\sqrt{\frac{\beta(t)}{2}}$. It is true as
\begin{equation*}
\begin{split}
& x_1(t) \leq 1-\sqrt{\frac{\beta(t)}{2}}\iff \beta(t)\leq \frac{(1-x_1(t))^2}{2}  \\
\Longrightarrow\; & \beta(t)< \frac{(1-x_1(t))^2}{2-x_1(t)}\iff \beta(t)< (1-x_1(t))^2-\beta(t)(1-x_1(t))\\
\Longrightarrow\; & \beta(t)< (1-x_1(t))^2-\beta(t)(1-x_1(t)) +\sum_{j=2}^r x_j(t)^2
\iff \beta(t)<\sum_{j=2}^r x_j(t) (1-x_j(t)-x_1(t))\\
\iff & \sum_{j=2}^r x_j(t) (1-x_j(t))> \sum_{j=2}^r x_j(t) x_1(t)+\beta(t) = (\sum_{j=2}^r x_j(t) +\beta(t))(1- \sum_{j=2}^r x_j(t)) > (1-x_1(t))x_1(t),
\end{split}
\end{equation*}
which implies~\eqref{eq:star_r_e3}. The asymptotic convergence of $x(t)$ to $\vect{e}_1$ 
is then established with the similar arguments in the proof of 
Theorem~\ref{thm:row-stochastic-star}~\rm{(ii)}.

Regarding fact~\rm{(iii)}, the existence and uniqueness of non-vertex equilibrium 
$x^*$ is established in the same way as in Theorem~\ref{thm:row-stochastic-general}~\rm{(i)}. 
$x^*$ satisfies~\eqref{eq:thm-e0} as well. The convergence property is similar to 
that of Theorem~\ref{thm:row-stochastic-general}~\rm{(ii)}. Specifically, $x_i(t)>0$ 
for all $t\geq n$ and $1\in\until{r}$. If we write $W(t)$ in the normal form as 
in~\eqref{eq:w-2g}, the statements (A1) and (A2) in the proof of 
Theorem~\ref{thm:row-stochastic-general}~\rm{(ii)} holds for $W_1(t)$ by the same arguments. 
That implies that  $\prod_{k=0}^t W_1(t)$ converge exponentially to a rank--$1$ matrix with 
positive identical rows, which is equal to $\vectorones[r] (x_{(1,\cdots, r)}^*)^\top$ and 
$x_{(1,\cdots, r)}^*$ is determined by~\eqref{eq:thm-e0}. Denote  
$\prod_{k=0}^t W(t)= \begin{bmatrix} P_1(t)& 0 \\ P_{21}(t)& P_{22}(t) \end{bmatrix}$. It is 
clear that $P_1(t)=\prod_{k=0}^t W_1(t)$, $P_{22}(t)=\prod_{k=0}^t W_{22}(t)$ and 
$P_{21}(t)=P_{21}(t-1)W_1(t)+P_{22}(t-1)W_{22}(t)$. As $P_{22}(t)$ converges exponentially to 
$\vectorzeros[]$, $P_{21}(t)$ then exponentially converges to $\vectorones[n-r] (x_{(1,\cdots, r)}^*)^\top$ 
following the previous statement that $\prod_{k=0}^t W_1(t)$ exponentially converges to 
$\vectorones[r] (x_{(1,\cdots, r)}^*)^\top$. Overall, $\prod_{k=0}^t W(t)$ converge 
exponentially to a rank--$1$ matrix with identical rows such that $x_{(1,\cdots, r)}$ converges 
exponentially to $x_{(1,\cdots, r)}^*>0$, and $x_{(r+1,\cdots, n)}$ converges exponentially to 
$\vectorzeros[n-r]$, for any non-vertex $x(0)$.
\end{IEEEproof}

\section{Proof of Theorem~\ref{thm:reducible-no-reachable}}
\label{pf:thm:reducible-no-reachable} 

\begin{IEEEproof}
Regarding the first part of fact~\rm{(i)}, the result is directly from the 
definition of the single-timescale DF model and has been proved in Theorem~\ref{thm:row-stochastic-star} and
Theorem~\ref{thm:reducible-reachable}: 
$x(0)=\vect{e}_i$ implies $x(t)=x^*=\vect{e}_i$ for all $t\geq 0$ and $i\in \until{n}$. 

Regarding fact~\rm{(i.1)}, as we discussed in the proof of 
Theorem~\ref{thm:reducible-reachable}, $\prod_{\tau=0}^t W_{MM}(\tau)$ 
converges exponentially to $\vectorzeros[m\times m]$ as $t$ goes to infinity, 
given $x_i(0)<1$ for all reducible node $i$. That implies that $x_{MM}(t)$ 
converges exponentially to $\vectorzeros[m]$ as $t$ goes to infinity. 

Regarding fact~\rm{(i.2)}, on an equilibrium $x^*$ in a sink $k$ with only two nodes, 
it shall satisfy from~\eqref{eq:sys_IF} that 
\begin{equation} \label{eq:reducible-no-reachable-e1}
x_{{kk}_1}^* (1 - x_{{kk}_1}^*)=x_{{kk}_2}^* (1 - x_{{kk}_2}^*).
\end{equation}
If $\zeta_k^*<1$, then the only solution to~\eqref{eq:reducible-no-reachable-e1} 
is  $x_{{kk}_1}^*=x_{{kk}_2}^*=\zeta_k^*/2$. If $\zeta_k^*=1$, then any pair 
$(\alpha,1-\alpha)^\top$ satisfies~\eqref{eq:reducible-no-reachable-e1} and hence, 
$x^*_{kk}=(\alpha,1-\alpha)^\top$ where $\alpha\in[0,1]$ depends upon the initial 
conditions and the topology of the network.

Regarding fact~\rm{(i.3)} and fact~\rm{(i.4)}, the proof 
is similar to the analysis of Theorem~\ref{thm:row-stochastic-general}~\rm{(i)}. 
For $C_{kk}$ irreducible and $\max \{c_{{kk}_i}\}<0.5$ we will show that there exists 
a unique $x_{kk}^*\in \interior{\simplex{n}}$ satisfying 
$x_{kk}^*=W_{kk} (x_{kk}^*)^\top x_{kk}^*+ W_{Mk}(x_{MM}^*)^\top x_{MM}^*=0$. 
As $x_{MM}^*=0$ from fact~\rm{(i.1)} above, the fix points shall satisfy
$x_{kk}^*-{x_{kk}^*}^2= C_{kk}^\top (x_{kk}^*-x_{kk}^*)$. As $x_{kk}^*$ shall 
be real valued and non-negative, given $\zeta_k^*>0$, $(x_{kk}^*-{x_{kk}^*}^2)$ 
is a scalar multiple of the unique positive left eigenvector of $C_{kk}$ associated 
with eigenvalue $1$.  
As $\max \{c_{{kk}_i}\}<0.5$ and $n_k\geq 3$, 
\begin{equation} \label{eq:reducible-no-reachable-e2}
x_{kk}^*-{x_{kk}^*}^2= \alpha_{kk}^* c_{kk}, \;\mbox{or equivalently,} 
\; x_{{kk}_i}^*=\alpha_{kk}^*\frac{c_{{kk}_i}}{1-{x_{{kk}_i}^*}}, 
\; \text{for all }i\in \until{n_k}, 
\end{equation}
where the scalar $\alpha_{kk}^*$ is such that 
$\vectorones[n_k]^\top x_{kk}^*= \zeta_k^*$, that is to say, 
\begin{equation*}
\alpha_{kk}^*=\frac{\zeta_k^*}{\sum_{j=1}^n c_{{kk}_j}/(1-x_{{kk}_j}^*)}.  
\end{equation*}
One may check that this $x_{kk}^*$ have the same form as the non-autocratic fixed 
point we obtained from the DF model~\cite{PJ-NEF-FB:14m}. Therefore, the uniqueness of $x_{kk}^*$ is 
directly from Theorem~3.6 in~\cite{PJ-NEF-FB:14m}. Moreover, the ordering of 
the elements of $x_{kk}^*$ is consistent with that of $c_{kk}$ 
following~\eqref{eq:reducible-no-reachable-e2}. 

Regarding fact~\rm{(ii)}, as $x_{kk}(t+1)=W_{kk}(t)^\top x_{kk}(t)+ W_{Mk}(t)^\top x_{MM}(t)$
with $W_{kk}(t)$ row stochastic and $W_{Mk}(t)^\top x_{MM}(t)\geq 0$, it is clear that
$\vectorones[n_k]^\top x_{kk}(t+1)\geq \vectorones[n_k]^\top W_{kk}(t)^\top x_{kk}(t)$.
That is $\zeta_{k}(t+1)\geq \zeta_{k}(t)$. For the second statement in fact~\rm{(ii)}, 
subject to Assumption 1) $\zeta_k(0)=0$ and Assumption 2) $x_i(0)=0$ for any reducible node 
$i$ such that there exists a directed path from $i$ to the sink $k$ in the network, we have 
$x_{kk}(1)=W_{kk}(0)^\top x_{kk}(0)+ W_{Mk}(0)^\top x_{MM}(0)=0$ as 
$W_{Mk}(0)^\top x_{MM}(0)>0$ contradicts thef second assumption above. Iteratively, we have
$x_{kk}(t+1)=W_{kk}(t)^\top x_{kk}(t)+ W_{Mk}(t)^\top x_{MM}(t)=0$ for all $t\geq 0$, where
the second term shall be equal to $\vectorzeros[m]$ for all the time as, otherwise it contradicts the 
second assumption.

Regarding fact~\rm{(iii)}, 
%first considering $x(0)=\vect{e}_i$ for 
%all $i\in \until{n}$, the self-weight configuration $x(t)$ of the dynamical 
%system~\eqref{eq:sys_IF} will remain on that vertex of the simplex for all times 
%as we have proved in fact~\rm{(i)}. 
%
%Second, we consider $x(0)\neq \vect{e}_i$ for any $i\in \until{n}$. We 
we will consider the convergence behaviors of self-weights in three different 
scenarios as described in facts~\rm{(i.1)}--~\rm{(i.3)}. 
%It is noted 
%that the convergence in each sink is independent from those of the rest of sinks 
%due to the connectivity of the network and the system setup~\eqref{eq:sys_IF}.

Scenario 1: The exponential convergence of the self-weights on 
reducible nodes has been clarified in fact~\rm{(i.1)}. 

Scenario 2: The convergence of the self-weights on a sink with only two nodes 
is similar to that described in Theorem~\ref{thm:reducible-reachable} fact~\rm{(i)} or fact~\rm{(iii)}. 
The difference is that all self-weights are accumulated on the two irreducible 
nodes in Theorem~\ref{thm:reducible-reachable} fact~\rm{(i)} 
but here $\zeta_k^*$ may be less than $1$ depending upon the initial condition 
and the topology of the network. If $\zeta_k^*=1$, then the convergence process 
here is exactly the same as Theorem~\ref{thm:reducible-reachable} 
fact~\rm{(i)}. If $\zeta_k^*<1$, then the self-weights in the two-node sink here 
exponentially converge to a unique $x_{kk}^*$. The analysis is 
similarly to that in Theorem~\ref{thm:reducible-reachable} 
fact~\rm{(iii)}. As these two nodes have the same eigenvector centrality score, 
the unique equilibrium is $(\zeta_k^*/2, \zeta_k^*/2)^\top$ here.

Scenario 3: The convergence of the self-weights on a sink with three or more 
nodes is almost the same as that described in 
Theorem~\ref{thm:reducible-reachable} fact~\rm{(iii)}. The only difference is 
that all self-weights are accumulated on the irreducible nodes as in 
Theorem~\ref{thm:reducible-reachable} fact~\rm{(iii)}
but here $\zeta_k^*$ may be less than $1$ depending upon the initial 
condition and the topology of the network.
If $\zeta_k^*=1$, then the analysis is the same to that of 
Theorem~\ref{thm:reducible-reachable} fact~\rm{(iii)} 
or that of Theorem~\ref{thm:row-stochastic-general} fact~\rm{(ii)}. 
If $\zeta_k^*<1$, then we have $\zeta_k(t)$ are upper bounded away from 
$1$ for all time $t$. As $\zeta_k(t)$ is non-decreasing, if $\zeta_k(t)>0 $ for 
$t=m$ (i.e., the max time for the social power migrating from a reducible node to the sink) then $\zeta_k(t)$ is uniformly bounded away from $0$ for all 
$t\geq m$, otherwise, if $\zeta_k(m)=0$ then $\zeta_k^*=0$. 
Given $\zeta_k(t)$ bounded away from $1$ and $0$ uniformly, 
first we have $x_{kk}(t)>0$ and each $x_{{kk}_i}$ is bounded way from $1$. 
Second, there exists a time $\tau$ such that any node in this sink has its 
self-weight $x_{{kk}_i}(t)$ lower bounded away from $0$ for all $t\geq \tau$. 
If it is not true, then by the irreducible property of $W_{kk}$ and the system 
definition~\eqref{eq:sys_IF}, all its connected nodes (i.e., all nodes in the 
sink) shall be sufficiently close to $0$ or $1$ for infinite time instances (see 
the similar argument~\eqref{eq:general-proof-e6} in the proof of 
Theorem~\ref{thm:row-stochastic-general}), that implies that $\zeta_k(t)$ is 
sufficiently close to $0$ or $1$ for infinite time instances, which is a 
contradiction. 
Third, the sum $\zeta_k(t)$ of the self-weights in this sink 
converges once all self-weights on reducible nodes exponentially converge 
to $0$, and the self-weight dynamics in the sink are independent from the 
dynamics occurred in other sinks. 
Finally, we can conclude that the exponential convergence of the product 
of $W_{kk}(t)$ based upon all results above. Consequently, $x_{kk}(t)$ 
converges exponentially as we have shown similarly in 
Theorem~\ref{thm:reducible-reachable} fact~\rm{(iii)}. 
\end{IEEEproof}

\end{document}